\input ./style/arxiv-general.cfg
\documentclass[aos,MSNbibl,nameyear,seceqn,dvips]{arximspdf}
\makeatletter
   \@ifpackageloaded{graphicx}{}{\usepackage{graphicx}}
\makeatother
\doi{10.1214/15-AOS1352}% Updated by VTEXPTS2LaTeX.exe, 09.07.2015
\volume{43}
\issue{6}
\pubyear{2015}
\firstpage{2565}
\lastpage{2587}
\docsubty{FLA}

\makeatletter
\def\mid{|}

\newproclaim{definition}{Definition}[section]
\newproclaim{example}{Example}[section]
\newtheorem{theorem}{Theorem}[section]
\newtheorem{lemma}{Lemma}[section]
\newtheorem{corollary}{Corollary}[section]

\newcommand{\betab}{\bolds{\beta}}
\newcommand{\mub}{\bolds{\mu}}
\newcommand{\etab}{\bolds{\eta}}
\newcommand{\deltab}{\bolds{\delta}}
\newcommand{\zetab}{\bolds{\zeta}}
\newcommand{\Jb}{\mathbf{J}}
\newcommand{\Vb}{\mathbf{V}}
\newcommand{\db}{\mathbf{d}}
\newcommand{\Cb}{\mathbf{C}}
\newcommand{\Ab}{\mathbf{A}}
\newcommand{\Bb}{\mathbf{B}}
\newcommand{\Xb}{\mathbf{X}}
\newcommand{\Ib}{\mathbf{I}}
\newcommand{\Mb}{\mathbf{M}}
\newcommand{\Wb}{\mathbf{W}}
\newcommand{\Ub}{\mathbf{U}}
\newcommand{\yb}{\mathbf{y}}
\newcommand{\Fb}{\mathbf{F}}
\newcommand{\Eb}{\mathbf{E}}
\newcommand{\Hb}{\mathbf{H}}
\newcommand{\Pb}{\mathbf{P}}
\newcommand{\zerob}{\mathbf{0}}
\newcommand{\epsb}{\bolds{\varepsilon}}
\def\red{}
\makeatother

\begin{document}
\begin{frontmatter}

%\dochead{}
\title{Optimal experimental designs for fMRI via circulant biased
weighing designs}
\runtitle{fMRI designs via circulant weighing designs}

\begin{aug}
% Corresponding author: Ming-Hung Kao - mkao3@asu.edu% Updated by
%VTEXPTS2LaTeX.exe, 09.07.2015 08:19
\author[A]{\fnms{Ching-Shui}~\snm{Cheng}\thanksref{M1,M3}\ead[label=e1]{cheng@stat.sinica.edu.tw}}
\and
\author[B]{\fnms{Ming-Hung}~\snm{Kao}\corref{}\thanksref{T1,M2}\ead[label=e2]{mkao3@asu.edu}}
\runauthor{C.-S. Cheng and M.-H. Kao}
\thankstext{T1}{Supported by NSF Grant DMS-13-52213.}
\affiliation{Academia Sinica\thanksmark{M1}, University of
California, Berkeley\thanksmark{M3}\\ and Arizona State
University\thanksmark{M2}}
%\author[A]{\fnms{}~\snm{}\corref{}\ead[label=e1]{}}%,
%\author[]{\fnms{}~\snm{}\ead[label=]{}}
% \and
%\author[]{\fnms{}~\snm{}\ead[label=]{}}
%\runauthor{}
%\affiliation{}
%\dedicated{}
\address[A]{Institute of Statistical Science\\
Academia Sinica\\
Taipei 11529\\
Taiwan\\
%Department of Statistics\\
%University of California, Berkeley\\
%Berkeley, CA 94720\\
\printead{e1}}

\address[B]{School of Mathematical and\\
\quad Statistical Sciences\\
Arizona State University\\
Tempe, Arizona 85287\\
USA\\
\printead{e2}}
%\address[A]{\\\printead{e1}}
%\address[]{\\\printead{}}
\end{aug}

% HISTORY:
%
\received{\smonth{12} \syear{2014}}% Updated by VTEXPTS2LaTeX.exe,
%09.07.2015 08:19
%
\revised{\smonth{6} \syear{2015}}% Updated by VTEXPTS2LaTeX.exe,
%09.07.2015 08:19

% ABSTRACT
%
\begin{abstract}
Functional magnetic resonance imaging (fMRI) technology is popularly
used in many fields for
studying how the brain reacts to mental stimuli. The identification of
optimal fMRI experimental designs is crucial for rendering
precise statistical inference on brain functions, but research on this
topic is very lacking. We
develop a general theory to guide the
selection of fMRI designs for estimating a hemodynamic response
function (HRF) that models the effect over time of the mental stimulus,
and for studying the comparison of two
HRFs. We provide a useful connection between fMRI designs
and circulant biased weighing designs, establish the statistical
optimality of some well-known fMRI designs and identify several new classes
of fMRI designs. Construction methods of high-quality fMRI designs are
also given.
\end{abstract}

% KEYWORDS
% Pirmas kwd is didziosios raides
%
\end{frontmatter}

%s1 #&#
\section{Introduction} \label{Sec1}
The present study concerns important issues on the design of
neuroimaging experiments where the pioneering functional magnetic
resonance imaging (fMRI) technology is employed to gain better
knowledge on how our brain reacts to mental stimuli. In such an fMRI
study, a sequence of tens or hundreds stimuli (e.g., images of
1.5-second flickering checkerboard) is presented to a human subject
while an fMRI scanner repeatedly scans the subject's brain to collect
data for making statistical inference about brain activity; see \citet
{Lazar2008bk51} for an overview of fMRI. The quality of such an
inference largely hinges on the amount of useful information contained
in the data, which in turn depends on the selected stimulus sequence
(i.e., fMRI design). The importance of identifying high-quality
experimental designs for fMRI and gaining insights into these designs
cannot be overemphasized.

In a seminal work, \citet{BuracasBoynton2002ar51} advocated the
use of the maximal length linear feedback shift
register sequences (or $m$-sequences) as fMRI designs for a precise
estimate of the hemodynamic response function (HRF);
the HRF models the effect over time of a brief stimulus on the relative
concentration of oxy- to deoxy-blood in the cerebral blood
vessels at a brain voxel (3D image unit), and is often studied for
gaining information about the underlying brain activity evoked by the stimulus.
The $m$-sequences have since then become very popular in practice. They are
also included as part of the ``good'' initial designs in the computer
algorithm of \citet{Kaoetal2009ar51} to facilitate the search of
multi-objective fMRI designs.
The computational results of \citet{BuracasBoynton2002ar51} and
\citet{Liu2004ar51} suggested
that the $m$-sequences can yield high statistical efficiencies in terms
of the $A$-optimality criterion; the $A$-criterion, which measures the
average variance of parameter estimates, is a design selection
criterion widely used in many fields including fMRI [e.g., \citet{Dale1999ar51, Fristonetal1999ar51}].
By focusing on the $D$-optimality criterion, \citet{Kao2014ar51}
proved the statistical optimality of
the binary $m$-sequences in estimating the HRF. As indicated there, the
binary $m$-sequences are special cases of the \textit{Hadamard
sequences} that can be generated from a certain type of Hadamard
matrices or difference sets (Section~\ref{sec2.3}); all these
designs are \mbox{$D$-optimal} in the sense of minimizing the generalized
variance of the HRF parameter estimates. While these designs are
expected to be $A$-optimal, there unfortunately is no theoretical proof
of this. One of our contributions here is to address this void. We also
identify some new classes of optimal fMRI designs for the estimation of
the HRF.

Another common study objective is on comparing HRFs of two stimulus
types (e.g., pictures of familiar vs. unfamiliar faces). Some
computational results on optimal fMRI designs for this study objective
have been reported in \citet{WagerNichols2003ar51}, \citet
{Liu2004ar51}, \citet{Kaoetal2008ar51}, \citet{Kaoetal2009ar51} and \citet{Mausetal2010ar51}. However,
theoretical work on providing insightful knowledge to guide the
selection of designs is scarce. In their pioneering papers, \citet
{LiuFrank2004ar51} and \citet{Mausetal2010ar51}
approximated the frequency of each stimulus type that an $A$- or
$D$-optimal fMRI design should possess. However, designs attaining the
optimal stimulus frequency can still be sub-optimal since the onset
times and presentation order of the stimuli play a vital role. Working
on this research line, \citet{Kao2014ar52} provided a sufficient
condition for fMRI designs
to be universally optimal in the sense of \citet
{Kiefer1975ar2021}, and proposed to construct optimal designs for
comparing two HRFs via an extended $m$-sequence (or de Bruijn
sequence), a Paley difference set or a circulant partial Hadamard
matrix. A major limitation of this recent contribution is that the
proposed designs exist only when the design length $N$ is a multiple of
$4$. New developments on identifying optimal fMRI designs for other
practical $N$ are called for.

We consider the two previously described design problems to target
optimal fMRI designs for the estimation of the HRF of a stimulus type
and for the comparison of two HRFs. Our main idea for tackling these
design issues is by formulating them into \textit{circulant biased
weighing design} problems. With this approach, we are able to prove
that the Hadamard sequences are optimal in terms of a large class of
optimality criteria that include both $A$- and $D$-criteria. This holds
as long as the design length $N (=4t -1$ for a positive integer $t$) of
such a design is sufficiently greater than the number of the HRF
parameters, $K$. For given $K$, a lower bound of $N$ for the design to
be both $A$- and $D$-optimal is also derived. This bound is easily
satisfied in typical fMRI experiments. In addition, we adapt and extend
previous results on (biased) weighing designs to identify some optimal
fMRI designs for estimating the HRF when $N = 4t$ and $N = 4t + 1$.
These results are further extended to cases where the study objective
is on the comparison of two HRFs. We note that the designs that we
present here exist in many design lengths for which optimal fMRI are
hitherto unidentified. These designs can be applied in practice or
serve as benchmarks to evaluate other designs; they help to enlarge the
library of high-quality fMRI designs.

The remainder of the paper is organized as follows. In Section~\ref{sec2}, we
provide relevant background information and present our main results on
optimal fMRI designs for estimating
the HRF. Our results on optimal fMRI designs for the comparison of two
HRFs are in Section~\ref{sec3}, and a conclusion is in Section~\ref{sec4}. Some proofs of
our results are presented in the \hyperref[app]{Appendix}.

%s2 #&#
\section{Designs for estimating the HRF} \label{sec2}
%s2.1 #&#
\subsection{Statistical model and design selection criteria} \label{sec2.1}
Consider an fMRI study where a mental stimulus such as a 1.5-second
flickering checkerboard
image [\citet{Boyntonetal1996ar51, Miezinetal2000ar51}] or a
painful heat stimulus [\citet{Worsleyetal2002ar51}]
is presented/applied to a subject at some of the $N$ time points in the
experiment. Let $y_n$ be the measurement of a brain voxel
collected by an fMRI scanner at the $n$th time point, $n=1,\ldots,N$.
We consider the following statistical model:
%
%e2.1 #&#
\begin{equation}
\qquad y_n = \gamma+ x_n h_1 + x_{n-1}
h_2 + \cdots+ x_{n-K+1} h_K +
\varepsilon_n,\qquad  n = 1, \ldots, N. \label{Model1}
\end{equation}
Here, $\gamma$ is a nuisance parameter, $h_1$ represents the unknown
height of the hemodynamic response function, HRF,
at the stimulus onset time point and $h_k$ is the HRF height at the
$(k-1)${th} time point following the stimulus onset.
The pre-specified integer $K$ is such that the HRF becomes negligible
after $K$ time points. The value of $x_{n-k+1}$ in
model (\ref{Model1}) is set to $1$ if $h_k$ contributes to $y_n$ and
$x_{n-k+1} = 0$ otherwise, and $\varepsilon_n$ is noise.

Our first design goal is to find an fMRI design, $\db= (d_1, \ldots,
d_N)^T$, that allows the most precise least-squares estimate of the HRF
parameter vector, $\mathbf{h}=(h_1, \ldots, h_K)^T$; here,
$d_n=1$ when a
stimulus appears at the $n$th time point and $d_n =0$ indicates no
stimulus presentation at that time point, $n=1,\ldots,N$. For simplicity,
we adopt the following assumptions from previous studies [\citet{Kao2013ar51} and
references therein]; see also Kao (\citeyear{Kao2014ar51,Kao2014ar52}) for discussions on these
assumptions. First, the last $K-1$ elements of $\db$ are also presented
in the pre-scanning period, that is, before the collection of $y_1$.
With this assumption, the value of $x_n$ in model (\ref{Model1}) is
$d_n$ for $n =1,\ldots, N$, and $x_{n} = d_{N+n}$ for $n \leq0$. In
addition, while additional nuisance terms may be included in the model
at the analysis stage to, say, allow for a trend/drift of $\yb= (y_1,
\ldots, y_N)^T$, we do not assume this extra complication when deriving
our analytical results on identifying optimal designs. We also consider
independent noise, but our results remain true when $\operatorname{cov}(\epsb) =
\alpha\Ib_N + \betab\mathbf{j}^T_N + \mathbf{j}_N\betab^T$, where $\epsb=(\varepsilon
_1, \ldots, \varepsilon_N)^T$, $\alpha$ is a constant, $\betab$ is a
constant vector, $\Ib_{N}$ is the identity matrix of order $N$, and
$\mathbf{j}
_N$ is the vector of $N$ ones; see also \citet{Kushner1997ar51}.
Other correlation structures of $\epsb$ such as an autoregressive
process may be considered, and is a focus of our future study. We now
rewrite model (\ref{Model1}) in the following matrix form:
%
%e2.2 #&#
\begin{equation}
\yb= \gamma\mathbf{j}_N + \Xb_d \mathbf{h}+ \epsb, \label{Eq:Model}
\end{equation}
where $\Xb_d = [\db, \Ub\db, \ldots, \Ub^{K-1}\db]$, and
%
%e2.3 #&#
\begin{equation}
\Ub= \left[ \matrix{ \zerob^T_{N-1}
& 1
\vspace*{2pt}\cr
\Ib_{N-1} & \zerob_{N-1}}\right].
\label{Eq:U}
\end{equation}
The information matrix for $\mathbf{h}$ is
$\Mb_b(\Xb_d) = \Xb_d^T(\Ib_N - N^{-1}\Jb_N)\Xb_d$, where $\Jb_N
= \mathbf{j}
_N\mathbf{j}^T_N$. We also let $\Mb(\Xb_d)= \Xb_d^T\Xb_d$.
Our target is at a $\db\in\mathcal{D}
= \{0,1\}^N$
that minimizes some real function $\Phi\{\Mb_b(\Xb_d)\}$ of $\Mb
_b(\Xb
_d)$. We consider the $A$-optimality criterion, $\Phi_A\{\Mb\} = \operatorname{tr}\{
\Mb
^{-1}\}/K$ for a positive definite $\Mb$, and $D$-optimality criterion,
$\Phi_D\{\Mb\} = |\Mb|^{-1/K}$. In addition, we adopt below some other
notions of optimality of designs and information matrices.
Specifically, the universal optimality described in Definition~\ref
{Def:UO} is due to \citet{Kiefer1975ar2021}. The type 1 criteria
of \citet{Cheng1978ar2021} with the version of \citet
{Cheng2014ar2021}, the $\Phi_p$-optimality criteria of \citet
{Kiefer1974ar2021} for $p \geq0$, and the $(M,S)$-optimality
[\citet{EcclestonHedayat1974ar2021}] are also considered.
Throughout this work, we set the criterion value to $+\infty$ for
designs with a singular information matrix.

%de2.1 #&#
\begin{definition}\label{Def:UO}
A design $\db$ is said to be universally optimal over a design class if
it minimizes $\Phi\{\Mb_b(\Xb_d)\}$ for all convex functions $\Phi$
such that
(i) $\Phi\{c\Mb\}$ is nonincreasing in $c>0$ and (ii) $\Phi(\Pb\Mb
\Pb
^T) = \phi(\Mb)$ for any $\Mb$ and any orthogonal matrix $\Pb$.
\end{definition}

%de2.2 #&#
\begin{definition}\label{Def:Type1} A design $\db$ is said to be
optimal over a design class with respect to all the type 1 criteria if
it minimizes $\Phi_{(f)}\{\Mb_b(\Xb_d)\} =  \sum_{i=1}^K f(\lambda
_i(\Mb
_b (\Xb_d)))$ for any real-valued function $f$ defined on $[0, \infty)$
such that
(i) $f$ is continuously differentiable in $(0, \infty)$ with $f' <0$,
$f'' > 0$, and $f''' <0$ and (ii) $\lim_{x \rightarrow0^+} f(x) = f(0)
= \infty$. Here $\lambda_i(\Mb_b(\Xb_d))$ is the $i$th greatest
eigenvalue of $\Mb_b(\Xb_d)$, $i=1,\ldots,K$.
\end{definition}

%de2.3 #&#
\begin{definition}\label{Def:Phip} A design $\db$ is said to be $\Phi
_p$-optimal over a design class for a given $p \geq0$ if
it minimizes
\[
\Phi_p\bigl\{\Mb_b(\Xb_d)\bigr\} = \cases{
\bigl|\Mb_b(\Xb_{d})\bigr|^{-1/K},
& \quad $\mbox{for $p=0$;}$
\vspace*{2pt}\cr
\bigl[ \operatorname{tr}\bigl\{\Mb^{-p}_b(\Xb_d)\bigr\}/K
\bigr]^{1/p}, &\quad $\mbox{for $p \in(0, \infty)$;}$
\vspace*{2pt}\cr
\lambda_1\bigl(\Mb^{-1}_b(
\Xb_d)\bigr), &\quad  $\mbox{when $p=\infty$,}$ }
\]
where $\lambda_i(\Mb_b(\Xb_d))$ is defined as in Definition~\ref{Def:Type1}.
\end{definition}

%de2.4 #&#
\begin{definition}\label{Def:MSopt} A matrix $\Mb^*$ is said to be
$(M,S)$-optimal over a class $\mathcal{M}$
of nonnegative definite matrices if (i) $\operatorname{tr}\{\Mb^*\} = \max_{M \in
\mathcal{M}} \operatorname{tr}\{\Mb\}$, and (ii) $\operatorname{tr}\{(\Mb^*)^2\} = \min_{M \in
\mathcal{M}_m} \operatorname{tr}\{\Mb^2\}$, where $\mathcal{M}_m \subset\mathcal{M}$
consists of all the matrices having the same trace as $\Mb^*$.
\end{definition}

Furthermore, we only consider optimality criteria $\Phi$ such that
%
%e2.4 #&#
\begin{equation}
\qquad\mbox{if}\quad \Phi(\Mb_1)\leq\Phi(\Mb_2),\qquad\mbox{then }
\Phi(c\Mb _1)\leq\Phi (c\Mb_2)\qquad\mbox{for all }c>0.
\label{Eq:c}
\end{equation}

%s2.2 #&#
\subsection{Circulant biased weighing designs} \label{sec2.2}
Our strategy for finding optimal fMRI designs is by taking advantage of
the link between these designs and circulant biased weighing designs.
A biased weighing design problem concerns the selection of a design for
efficient estimation of the weights of $K$ objects in $N$ weighings on
a spring/chemical balance
that has an unknown systematic bias. A~spring balance weighing design
(SBWD) is specified by
a $\Wb\in\{0, 1\}^{N\times K}$, where the $(n,k)${th} element of
$\Wb
$ indicates that the $k$th object is placed on the balance~($1$), or
absent ($0$) in the $n$th weighing. Such a design is called \textit
{circulant} if $\Wb$ is a circulant matrix. The information matrix
$\Mb
_b(\Wb)$ for the $K$ weights is equal to $\Wb^T(\Ib_N - N^{-1}\Jb
_N)\Wb
$. For each fMRI design $\db$, the matrix $\Xb_d$ clearly defines a
circulant SBWD. Thus, the fMRI design issue formulated earlier is a
sub-problem of the optimal SBWD problem: selecting an optimal design
among circulant SBWDs.

A chemical balance weighing design (CBWD) is specified by a $\bar{\Wb}
\in\{-1, 0,  1\}^{N\times K}$, where the $(n,k)${th} element of $\bar
{\Wb}$ indicates that the $k$th object is placed on the left pan
($-1$), right pan ($+1$), or absent ($0$) in the $n$th weighing.
Each SBWD matrix $\Wb$ can be transformed into a CBWD matrix $\tilde
{\Wb
}$ via $\tilde{\Wb}=\pm(\Jb_{N,K}-2\Wb)$, where $\Jb_{N,K} =
\mathbf{j}_N\mathbf{j}
_K^T$; that is, 0 and~1 are replaced by 1 and $-1$, or $-1$ and~1,
respectively. Given an fMRI design $\db\in\mathcal{D}$, let $\tilde
{\db
}=\pm(\mathbf{j}_N-2\db)$ and $\Xb_{\tilde{d}} = [\tilde
{\db}, \Ub\tilde{\db},
\ldots, \Ub^{K-1}\tilde{\db}]$, where $\Ub$ is defined as in (\ref
{Eq:U}); then, $\tilde{\db}\in\tilde{\mathcal{D}}=\{-1,1\}^N$, and
$\Xb
_{\tilde{d}}$ is a circulant CBWD matrix. Specifically, if we write
$\Xb
_d$ as $\Wb$, then $\Xb_{\tilde{d}}=\tilde{\Wb}$.

\citet{Cheng2014ar2021} showed that
%
%e2.5 #&#
\begin{equation}
\Mb_b(\Wb)=\tfrac{1}{4}\Mb_b(\tilde{\Wb})\qquad\mbox{for all }\Wb\in \{0, 1\} ^{N\times K}. \label{Eq:information}
\end{equation}
We thus have the following result.

%le2.1 #&#
\begin{lemma} \label{Lemma:biasedWD}
For any $\Phi$ satisfying (\ref{Eq:c}) and any $\Wb_1,\Wb_2\in\{0,
1\}
^{N\times K}$, $\Phi(\Mb_b(\Wb_1))\leq\Phi(\Mb_b(\Wb_2))$ if and only
if $\Phi(\Mb_b(\tilde{\Wb}_1))\leq\Phi(\Mb_b(\tilde{\Wb}_2))$.
Therefore,
an fMRI design $\db^* \in\mathcal{D}$ is $\Phi$-optimal for estimating
the HRF if and only if
\[
\Phi\bigl\{\Mb_b(\Xb_{\tilde{d}^*})\bigr\} = \min
_{\tilde{d} \in\tilde
{\mathcal
{D}}} \Phi\bigl\{\Mb_b(\Xb_{\tilde{d}})\bigr\}\qquad
\mbox{where } \tilde{\db }^* = \pm \bigl(\mathbf{j}_N - 2
\db^*\bigr).
\]
\end{lemma}

Lemma~\ref{Lemma:biasedWD} reduces the problem of finding optimal fMRI
designs for estimating an HRF to that of identifying optimal biased
circulant CBWDs without zero entries. In Section~\ref{sec3}, we
establish the connection of optimal designs for comparing two HRFs to
that of optimal biased circulant CBWDs allowing zero entries.

%s2.3 #&#
\subsection{Main results} \label{sec2.3}
Following Lemma~\ref{Lemma:biasedWD}, we tackle our first fMRI design
issue by working on circulant CBWDs ($\Xb_{\tilde{d}}$ for $\tilde
{\db}
\in\tilde{\mathcal{D}}$) that contain no zero. We have the following
result for such CBWDs.

%le2.2 #&#
\begin{lemma} \label{Lemma_mod3}
For $\tilde{\db} \in\tilde{\mathcal{D}}$, $\Mb(\Xb_{\tilde{d}})
= \Xb
^T_{\tilde{d}}\Xb_{\tilde{d}}$ has
diagonal elements equal to $N$ and off-diagonal elements congruent to
$N$ modulo 4. When $N$ is odd,
$\Mb_b(\Xb_{\tilde{d}}) \preceq\Mb(\Xb_{\tilde{d}}) - N^{-1}\Jb
_K$ and
the equality holds if $\mathbf{j}_N^T\tilde{\db} = \pm1$;
here, $\preceq$ is
the L\"owner ordering, that is, $\Mb_1 \preceq\Mb_2 $ if $\Mb_2 -\Mb
_1$ is
nonnegative definite.
\end{lemma}

\begin{pf}
All the diagonal elements of $\Mb(\Xb_{\tilde{d}})$ are clearly
$\tilde
{\db}^T\tilde{\db} = N$. In addition, for $q, r \in\{-1, +1\}$, let
$n^{(rq)}_k$ be the number of times $(\tilde{d}_{n-k}, \tilde
{d}_n) = (q,r)$, where $\tilde{d}_{n}$ is the $n$th element of
$\tilde{\db}$ for $n=1,\ldots, N$, and $\tilde{d}_n$ is set to
$\tilde
{d}_{N+n}$ when $n \leq0$. We have, for any $i \neq j$ and $k=|i-j|$,
the $(i,j)${th} element of $\Mb(\Xb_{\tilde{d}})$ is $n^{(++)}_k
+ n^{(--)}_k - (n^{(+-)}_k + n^{(-+)}_k) = N -
4n^{(-+)}_k$, and is thus congruent to $N$ modulo 4. Note that
the above equality is a consequence of the fact that $\Xb_{\tilde{d}}$
is a circulant matrix. Moreover, $\Mb_b(\Xb_{\tilde{d}}) = \Mb(\Xb
_{\tilde{d}}) - a^2 N^{-1}\Jb_K$, where $a = \mathbf{j}_N^T\tilde{\db}$ with
$a^2 \geq1$ if $N$ is odd. Thus, $\Mb(\Xb_{\tilde{d}}) - N^{-1}\Jb
_K -
\Mb_b(\Xb_{\tilde{d}}) = (a^2 - 1)N^{-1}\Jb_K$, and our claim follows.
\end{pf}

We now provide some results for obtaining optimal circulant biased
CBWDs with no zero, and hence, optimal fMRI designs for estimating the HRF.
For cases with $N=4t-1\ (\geq4$), the following lemma due to \citet
{Cheng1992ar2021} is useful.

%le2.3 #&#
\begin{lemma} \label{Lemma_C1992}
Let $\mathcal{M}^N$ be a set of $K$-by-$K$ symmetric and nonnegative
definite matrices,
$\mathcal{M}^N_m \subset\mathcal{M}^N$ be the set of matrices that
have the maximum trace over $\mathcal{M}^N$, and $\mathcal{M}^N_{ms}$
be the set of $\Mb$ that minimize $\operatorname{tr}(\Mb^2)$ over $\mathcal{M}^N_m$.
Suppose $A_N = \max_{M \in\mathcal{M}^N} \operatorname{tr}\{\Mb\}$ and $B_N = \min_{M
\in\mathcal{M}^N_m} \operatorname{tr}\{\Mb^2\}$ are such that \textup{(a)} $\lim_{N\rightarrow
\infty} A_N = \infty$, and \textup{(b)} for some $L >0$, $|B_N - K^{-1}A^2_N|
\leq L$ for all $N$. In addition, let $\lambda_i(\Mb)$ be as in
Definition~\ref{Def:Type1}, and
$\Phi_{(g)}\{\Mb\} = \sum_{i=1}^K g(\lambda_i(\Mb))$ for a
real-function $g$ satisfying the following two conditions: \textup{(i)} $g$ is
thrice continuously differentiable in a neighborhood of 1 with $g'(1) <
0$ and $g''(1) > 0$ and \textup{(ii)} for any $c > 0$, there are constants
$\alpha(c) > 0$ and $\beta(c)$ such that $g(cx) = \alpha(c) g(x) +
\beta
(c)$ for all $x$. Then there exists an $N_0(K, g)$ such that whenever
$N \geq N_0(K, g)$, for any $\Mb^*\in\mathcal{M}^N_{ms}$, we have
$\Phi
_{(g)}\{\Mb^*\} < \Phi_{(g)}\{\Mb\}$ for all $\Mb\notin\mathcal
{M}^N_{ms}$.
\end{lemma}

We note that when $g(x) = -\log x$, $\Phi_{(g)}$ is equivalent to the
$D$-criterion, or
equivalently, the $\Phi_p$-criterion in Definition~\ref{Def:Phip} with
$p=0$. This and the other $\Phi_p$-criteria satisfy conditions (i) and
(ii) in
Lemma~\ref{Lemma_C1992}; see also \citet{Cheng1992ar2021}. We thus
have the following result on $\Phi_p$-optimal circulant biased CBWDs and
$\Phi_p$-optimal fMRI designs for $N=4t -1$.

%th2.1 #&#
\begin{theorem} \label{thmm1}
Let $N = 4t - 1$, $p_0 > 0$, and $\tilde{\db}^* \in\tilde{\mathcal
{D}}$ be such that\break
$\Mb_b(\Xb_{\tilde{d}^*}) = (N+1)[\Ib_K - N^{-1}\Jb_K]$. Then there
exists an $N_0(K, p_0)$ such that, whenever
$N \geq N_0(K, p_0)$, $\tilde{\db}^*$ is $\Phi_p$-optimal over
$\tilde
{\mathcal{D}}$, and $\db^* = (\mathbf{j}_N \pm\tilde{\db
}^*)/2$ is $\Phi
_p$-optimal over $\mathcal{D}$
for any $p \in[0, p_0]$.
\end{theorem}

\begin{pf}
We first work on $\Mb(\Xb_{\tilde{d}}) - N^{-1}\Jb_K$ for $\tilde
{\db}
\in\tilde{\mathcal{D}}$. Following Lemma \ref{Lemma_mod3}, the
diagonal elements of $\Mb(\Xb_{\tilde{d}}) - N^{-1}\Jb_K$ are $N_b =
N-N^{-1}$, and the $(i,j)${th} element of this matrix is
$(c_{i,j} - N^{-1})$ with $c_{i,j} = 3 (\mathrm{mod}\ 4)$ for $i\neq j$. Thus,
$\operatorname{tr}\{[\Mb(\Xb_{\tilde{d}}) - N^{-1}\Jb_K]^2\}$ is minimized when
$c_{i,j} = -1$ for all $i \neq j$. This implies the $(M,S)$-optimality
of $\Mb(\Xb_{\tilde{d}^*}) - N^{-1}\Jb_K$ over $\tilde{\mathcal
{M}}_b =
\{\Mb(\Xb_{\tilde{d}}) - N^{-1}\Jb_K \mid\tilde{\db} \in\tilde
{\mathcal{D}} \}$. In addition, it can be seen that $A_N = KN_b$, and
$B_N = KN_b^2 +(1+N^{-1})^2K(K-1)$, where $A_N$ and $B_N$ are defined
as in Lemma~\ref{Lemma_C1992}. Therefore, $\lim_{N\rightarrow\infty}
A_N = \infty$, and $|B_N - K^{-1}A^2_N| = |KN_b^2 +(1+N^{-1})^2K(K-1) -
KN^2_b| = (1 + N^{-1})^2K(K-1)$ is bounded above by a positive number
for all $N$. Following Lemmas \ref{Lemma_mod3}, and \ref{Lemma_C1992},
we then have, when $N \geq N_0(K, p_0)$ for some $N_0(K, p_0)$,
\begin{eqnarray*}
\Phi_{p_0}\bigl\{\Mb_b(\Xb_{\tilde{d}^*})\bigr\} &=&
\Phi_{p_0}\bigl\{\Mb(\Xb _{\tilde
{d}^*}) - N^{-1}
\Jb_K\bigr\}
\\
& \leq&\Phi_{p_0}\bigl\{\Mb(\Xb_{\tilde{d}}) - N^{-1}
\Jb_K\bigr\} \leq\Phi_{p_0}\bigl\{\Mb_b(
\Xb_{\tilde{d}})\bigr\}
\end{eqnarray*}
for any $\tilde{\db} \in\tilde{\mathcal{D}}$; here $\Phi_{p_0}$ is
defined as in Definition~\ref{Def:Phip}. The $\Phi_p$-optimality of
$\tilde{\db}^*$ over $\tilde{\mathcal{D}}$ for $p \in[0, {p_0}]$ then
follows from Corollary~3.3 of \citet{Cheng1987ar2021} and the fact
that $\Mb_b(\Xb_{\tilde{d}^*})$ has two eigenvalues, with the smaller
one having multiplicity $1$. Moreover, with Lemma~\ref{Lemma:biasedWD},
we obtain the $\Phi_p$-optimality of $\db^*$ over $\mathcal{D}$.
\end{pf}

In Theorem~\ref{4546465545485}, we provide an $N_0(K, 1)$ for a design to be
$\Phi_1$-optimal (i.e., $A$-optimal).
Our approach for deriving this bound for $N$ is analogous to that of
\citet{GalilKiefer1980ar2021}, and
\citet{SatheShenoy1989ar2021}. A proof is provided in the \hyperref[app]{Appendix}.

%th2.2 #&#
\begin{theorem}\label{4546465545485}
Consider the same conditions as in Theorem~\ref{thmm1}. Let $N_0(K,  1)$
be the greatest real root of the cubic function
$c(x) = 2x^3 + (10-7K)x^2 + 2(2K-5)(K-1)x + 4K^2 - 7K$. If $K \geq4$
and $N \geq N_0(K, 1)$, then $\tilde{\db}^*$
is $\Phi_p$-optimal over $\tilde{\mathcal{D}}$, and $\db^*$ is
$\Phi
_p$-optimal over $\mathcal{D}$ for $0\leq p \leq1$.
\end{theorem}

Recently, \citet{Kao2014ar51} studied the efficiency of Hadamard
sequences, $\db_H$, in estimating $\mathbf{h}$ of model
(\ref{Eq:Model}).
A Hadamard sequence is a binary sequence constructed from a normalized
Hadamard matrix $\Hb\in\{-1,1\}^{(N+1) \times(N+1)}$ that contains a
circulant core $\tilde{\Hb}$. Such an $\Hb$ is such that $\Hb\Hb^T =
(N+1)\Ib_{N+1}$, the elements of its first row and column are all $1$,
and the bottom-right $N$-by-$N$ sub-matrix $\tilde{\Hb}$ is a circulant
matrix. These Hadamard matrices are known to exist when $N$ is a prime,
a product of twin primes, or $2^r - 1$ for an integer $r > 1$. They can
be easily generated by, for example, the Paley, Singer, or twin prime
power difference sets [\citet{GolombGong2005bk101,
Horadam2007bk101}]. Any column of the circulant core $\tilde{\Hb}$
is a vertex, $\tilde{\db}_H$, of the hypercube $\tilde{\mathcal{D}}
= \{
-1, 1\}^N$, and $\db_H = (\mathbf{j}_N -\tilde{\db}_H)/2$
forms a Hadamard
sequence. The popularly used binary $m$-sequences [\citet
{BuracasBoynton2002ar51}] can be constructed by the same method
when $N = 2^r - 1$, and are thus special cases of $\db_H$. The $\db_H$
has design length $N = 4t - 1$ for some integer $t$, and our results
can be applied to establish the $A$- and $D$-optimality of these
designs as stated in the following corollary.

%co2.1 #&#
\begin{corollary} \label{cor1}
A Hadamard sequences $\db_H$ of length $N$ is $A$- and $D$-optimal for
estimating the HRF if $K \geq4$ and $N \geq N_0(K, 1)$.
Here, $N_0(K, 1)$ is defined as in Theorem~\ref{4546465545485}.
\end{corollary}

\begin{pf} For $\db_H$, let $\tilde{\db}_H = \mathbf{j}_N- 2\db_H$. Then it
can be seen from the
construction of $\db_H$ that $\Xb_{\tilde{d}_H}$ is a ciruclant matrix
consisting of $K$ distinct columns of the circulant core of a
normalized Hadamard matrix. Consequently, $\Mb_b(\Xb_{\tilde{d}_H}) =
(N+1)[\Ib_K - N^{-1}\Jb_K]$.
Our claim then follows from Lemma~\ref{Lemma:biasedWD} and Theorem~\ref{4546465545485}.
\end{pf}

Our results so far are for cases with $N = 4t - 1$. For $N = 4t$, if
there exists a $\tilde{\db}$ with $\Mb_b(\Xb_{\tilde{d}}) = N\Ib_K$,
then $\tilde{\db}$ is universally optimal over $\tilde{\mathcal{D}}$,
and the corresponding $\db= (\mathbf{j}_N \pm\tilde{\db
})/2$ is universally
optimal in estimating the HRF over
all fMRI\vspace*{1pt} designs. This fact follows directly from Proposition~1$'$ of
\citet{Kiefer1975ar2021}. We note that $\tilde{\db}$ is
universally optimal whenever the columns of $\Xb_{\tilde{d}}$
are pairwise orthogonal, and are all orthogonal to $\mathbf{j}_N$. The
transpose of such a matrix $\Xb_{\tilde{d}}$ is called a \textit
{circulant partial Hadamard matrix} by \citet
{Craigenetal2013ar101}. Clearly, the corresponding $\Xb_{d}$ with
$\db= (\mathbf{j}_N \pm\tilde{\db})/2$ forms a
two-symbol, $N$-run,
$K$-factor \textit{circulant orthogonal array (OA)} whose strength is
$\geq2$; see \citet{Hedayatetal1999bk101} for an overview of
OAs. A circulant partial Hadamard matrix, and thus a circulant OA, can
be obtained by a computer search [\citet{Linetal2014ar51,
Lowetal2005ar101}]. Here, we provide a systematic method for
constructing a universally optimal $\db$.

%th2.3 #&#
\begin{theorem}\label{thmm3}
Let $\db_{1,g, H} \in\mathcal{D}$ be obtained by inserting a $0$ to a
run of $g$ $0$'s in a Hadamard sequence $\db_H$.
If $K \leq g + 1$, then $\db_{1,g,H}$ is universally optimal for
estimating $\mathbf{h}$ of model (\ref{Eq:Model}).
\end{theorem}

\begin{pf} Without loss of generality, we assume that a run of $g$
$0$'s appears in the tail of $\db_H$, and $\db_{1,g, H}$ is obtained by
adding a $0$ to this run of $0$'s. Suppose $K \leq q + 1$, and $\tilde
{\db}_{1,g, H} =\mathbf{j}_N - 2\db_{1,g, H}$. It can be
seen that $\Xb_{\tilde
{d}_{1,g, H}}$ is an $N$-by-$K$ circulant
%matrix
orthogonal array
whose columns are some $K$ distinct columns of a Hadamard matrix.
%, or more precisely, that of a cyclic shift of a normalized Hadamard
%matrix of order $N (=4t)$.
\end{pf}

For $N = 4t + 1$, Theorem~4.1 of \citet{Cheng2014ar2021} provides
a guidance on the selection of $\Phi_{(f)}$-optimal biased CBWDs for
any type 1 criterion $\Phi_{(f)}$.
We describe this result in Lemma~\ref{lemma5} with our notation. It is
interesting to note that, under our setting, a simple alternative proof
of Lemma~\ref{lemma5} can be achieved by
utilizing Theorem~2.1 of \citet{Cheng1980ar2021} that is slightly
rephrased in Lemma~\ref{C1980} below.

%le2.4 #&#
\begin{lemma}\label{C1980}
Let $\Mb^*$ be a symmetric matrix with eigenvalues $\lambda_1(\Mb^*) >
\lambda_2(\Mb^*) = \lambda_3(\Mb^*) = \cdots= \lambda_{K}(\Mb^*)
> 0$
and $\mathcal{M}$ be a set of nonnegative definite symmetric matrices.
If the following
conditions are satisfied, then\break $\Phi_{(f)}\{\Mb^*\} \leq\Phi_{(f)}\{
\Mb
\}$ for any $\Mb\in\mathcal{M}$ and any type 1 criterion $\Phi_{(f)}$:
\begin{longlist}[(a)]
\item[(a)] $\operatorname{tr}\{\Mb^*\} \geq \operatorname{tr}\{\Mb\}$ for any $\Mb\in\mathcal{M}$;
% \item[(b)] $\operatorname{tr}\{(\Mb^*)^2\} < (\operatorname{tr}\{\Mb^*\})^2/(K-1)$;
%
\item[(b)] for any $\Mb\in\mathcal{M}$, $\operatorname{tr}\{\Mb^*\} - \sqrt
{[K/(K-1)][\operatorname{tr}\{(\Mb^*)^2\} - (\operatorname{tr}\{\Mb^*\})^2/K]} \geq$ $\operatorname{tr}\{\Mb\} -
\sqrt{[K/(K-1)][\operatorname{tr}\{\Mb^2\} - (\operatorname{tr}\{\Mb\})^2/K]}$.
\end{longlist}
\end{lemma}

Note that since the condition $\lim_{x \rightarrow0^+} f(x) = f(0) =
\infty$ is required in Definition~\ref{Def:Type1}, there is no need to
verify (2.2) in Theorem~2.1 of \citet{Cheng1980ar2021}; see
Theorem~2.3 of \citet{Cheng1978ar2021}.

%t1 #&#
\begin{table}[b]
\caption{A Hadamard sequence $\db_H$, and a $\db_{1,g,H}$ for
estimating $\mathbf{h}$ with $K \leq9$}\label{tab1}
\begin{tabular*}{\textwidth}{@{\extracolsep{\fill}}lcl@{}}
\hline
% after \ \hline or \cline{col1-col2} \cline{col3-col4}...
& $\bolds{N}$ & \multicolumn{1}{c@{}}{\textbf{Design}} \\
\hline
$\db_H$ & 151 & 1 0 0 1 0 0 1 1 0 0 0 0 1 1 1 1 0 0 0 0 0 0 0 1 1 0 1 1 1 0 1 0 0 1 0 1\\
&&                0 0 0 0 0 1 0 0 0 0 1 0 1 0 0 1 1 1 1 0 1 1 0 0 1 1 0 1 0 1 1 1 0 0 1 1\\
&&                0 1 0 1 0 1 0 1 0 0 1 1 0 0 0 1 0 1 0 0 1 1 0 0 1 0 0 0 0 1 1 0 1 0 1 1\\
&&                1 1 0 1 1 1 1 1 0 1 0 1 1 0 1 0 0 0 1 0 0 1 1 1 1 1 1 1 0 0 0 0 1 1 1 1\\
&&                0 0 1 1 0 1 1 \\[3pt]
$\db_{1,g,H}$ & 132 & 1 0 1 0 0 0 1 0 1 0 1 0 0 0 1 0 0 1 1 1 0 0 1 1 1 0 1 0 0 1 1 1 1 0 0 0\\
&&                      0 1 0 0 1 0 1 0 0 0 0 1 0 0 1 1 0 0 1 0 1 1 0 0 0 0 0 0 0 0 0 1 1 1 1 1\\
&&                      1 1 1 0 0 1 0 1 1 0 0 1 1 0 1 1 1 1 0 1 0 1 1 0 1 1 1 1 0 0 0 0 1 1 0 1\\
&&                      0 0 0 1 1 0 0 0 1 1 0 1 1 1 0 1 0 1 0 1 1 1 0 1 \\
\hline
\end{tabular*}
\end{table}

%le2.5 #&#
\begin{lemma}\label{lemma5}
Let $N = 4t + 1$, and $\tilde{\db}^* \in\tilde{\mathcal{D}}$ be
such that
$\Mb_b(\Xb_{\tilde{d}^*}) = (N-1)[\Ib_K + N^{-1}\Jb_K]$.
Then $\tilde{\db}^*$ is optimal over $\tilde{\mathcal{D}}$, and
$\db^*
= (\mathbf{j}_N \pm\tilde{\db}^*)/2$ is optimal for
estimating the HRF in
terms of any type 1 criterion.
\end{lemma}

\begin{pf}
From Lemma~\ref{Lemma_mod3}, we have that the diagonal elements of\break
$\Mb
(\Xb_{\tilde{d}}) = \Xb^T_{\tilde{d}}\Xb_{\tilde{d}}$ are $N$,
and the
off-diagonal elements are congruent to $1$ modulo $4$.
In addition, $\Mb_b(\Xb_{\tilde{d}}) \preceq\Mb(\Xb_{\tilde{d}}) -
N^{-1}\Jb_K$, and the equality holds when $\mathbf{j}^T_N\tilde{\db} = \pm1$.
It can then be easily seen that $\Mb(\Xb_{\tilde{d}^*}) - N^{-1}\Jb_K$
is $(M,S)$-optimal over $\tilde{\mathcal{M}}_b = \{\Mb(\Xb_{\tilde{d}})
- N^{-1}\Jb_K \mid\tilde{\db} \in\tilde{\mathcal{D}}\}$, and
conditions (a) and (b) in Lemma~\ref{C1980} are satisfied if we replace
$\Mb^*$, $\Mb$ and $\mathcal{M}$ there by $\Mb(\Xb_{\tilde{d}^*}) -
N^{-1}\Jb_K$, $\Mb(\Xb_{\tilde{d}}) - N^{-1}\Jb_K$, and $\tilde
{\mathcal{M}}_b$,
respectively. Consequently, for any type 1 criterion $\Phi_{(f)}$, we have
$\Phi_{(f)}\{\Mb_b(\Xb_{\tilde{d}^*})\} = \Phi_{(f)}\{\Mb(\Xb
_{\tilde
{d}^*}) - N^{-1}\Jb_K\} \leq\Phi_{(f)}\{\Mb(\Xb_{\tilde{d}}) -
N^{-1}\Jb_K\}
\leq\Phi_{(f)}\{\Mb_b(\Xb_{\tilde{d}})\}$. The optimality of
$\tilde
{\db}^*$ thus follows. By (\ref{Eq:information}), this argument also
applies to $\db^*$.
\end{pf}

We now provide a systematic method for constructing optimal fMRI
designs for cases with $N = 4t + 1$, followed by an example
on an application of our results in this section.

%th2.4 #&#
\begin{theorem}\label{thmm4}
Let $\db_{2,g, H} \in\mathcal{D}$ be obtained by inserting two $0$'s
to a run of $g$ $0$'s in a Hadamard sequence $\db_H$.
If $K \leq g + 1$, then $\db_{2,g, H}$ is optimal for estimating
$\mathbf{h}$
of model (\ref{Eq:Model}) for all type 1 criteria.
\end{theorem}

\begin{pf} With a similar argument as in the proof of Theorem~\ref
{thmm3}, we have that, when $K \leq g+1$ and $\tilde{\db}_{2,g, H} =
\mathbf{j}
_N - 2\db_{2,g, H}$, $\Mb_b(\Xb_{\tilde{d}_{2,g, H}}) = (N-1)[\Ib
_K +
N^{-1}\Jb_K]$.
Our claim then follows from Lemma~\ref{lemma5}.
\end{pf}

%ex2.1 #&#
\begin{example}\label{examp1}
Consider an experiment where the stimulus can possibly occur every $4$
seconds. Then $N = 4(38) - 1 = 151$ corresponds
to a 10-minute experiment and $K = 9$ corresponds to a 32-second HRF. A
$\db_H =(d_{1,H},\ldots, d_{N,H})^T$ can be obtained by a Paley difference
set [\citet{Paley1933ar101}]. This is to set $d_{n,H}=0$ if $(n-1)
\in\{x^2\ (\mathrm{mod}\ N) \mid x = 1, \ldots, (N-1)/2\}$ and $d_{n,H}=1$,
otherwise. The obtained design $\db_H$ is presented in Table~\ref
{tab1}. It\vadjust{\goodbreak} is both $A$- and $D$-optimal for estimating $\mathbf{h}$ of model
(\ref{Eq:Model}) since $N > N_0(K, 1) = 21.34$.

Following Theorem~\ref{thmm3}, we may insert a $0$ to the longest run
of $0$'s in the $\db_H$ presented in Table~\ref{tab1} to yield a
universally optimal design. The resulting design can accommodate a $K
\leq8$. For $K = 9$, we obtain a universally optimal $\db_{1,g,H}$ by
extending a $\db_H$ of length $N = 131$. This $\db_{1,g,H}$ is
presented in Table~\ref{tab1}. We also obtain a $\db_{2,g,H}$ by
inserting another $0$ into the longest run of $0$'s in $\db_{1,g,H}$.
Following Theorem~\ref{thmm4}, this $\db_{2,g,H}$ is optimal for any
type 1 criterion in estimating $\mathbf{h}$ with $K \leq9$.
\end{example}

It is noteworthy that, by replacing $0$ and $1$ with $1$ and $2$,
respectively, the $\db_{1,g,H}$ in Table~\ref{tab1} is equivalent to
the design of the same design length in Table~3.1 of \citet
{Kao2014ar52}. The use of such a design whose elements are $1$ or
$2$ is discussed in the next section.

%% Comparing HRF
%s3 #&#
\section{Designs for contrasts between HRFs} \label{sec3}
We now consider optimal fMRI experimental designs for studies where the
objective is on comparing HRFs of two stimulus types.
For this situation, \citet{Kao2014ar52} presented some optimal
designs for $N=4t$ by considering the
following extension of model (\ref{Eq:Model}):
%
%e3.1 #&#
\begin{equation}
\yb= \gamma\mathbf{j}_N + \Xb_{u, 1} \mathbf{h}_1 + \Xb_{u, 2} \mathbf{h}_2
+ \epsb, \label{Eq:Model2}
\end{equation}
where $\mathbf{h}_q = (h_{q1}, \ldots, h_{qK})^T$ is the
vector of the $K$
unknown HRF heights of the $q${th}-type stimulus, $\Xb_{u,q}$
is the 0-1 design matrix obtained from the selected fMRI design $\mathbf{u}
=(u_1, \ldots, u_N)^T$ with $u_n \in\{0, 1, 2\}$, $q = 1, 2$,
and the remaining terms are as in (\ref{Eq:Model}). Specifically, $u_n
= q > 0$ indicates that
a stimulus of the $q${th} type appears at the $n$th time point, and
$u_n = 0$ if no
stimulus is present. In addition, for $q = 1, 2$, $\Xb_{u,q} =
[\deltab
_q, \Ub\deltab_q, \ldots, \Ub^{K-1}\deltab_q]$, where $\Ub$ is defined
in (\ref{Eq:U}), and the $n$th element of $\deltab_q$ is $1$ if $u_n
= q$, and is $0$ otherwise. The main interest lies in $\zetab= \mathbf{h}_1 -
\mathbf{h}_2$, and we may rewrite model (\ref{Eq:Model2}) as
%
%e3.2 #&#
\begin{equation}
\yb= \gamma\mathbf{j}_N + \Eb_{u}\etab+
\Fb_{u} \zetab + \epsb, \label{Eq:Model3}
\end{equation}
where $\Eb_{u} = (\Xb_{u,1} + \Xb_{u,2})/2$, $\etab=\mathbf{h}_1 + \mathbf{h}_2$, $\Fb
_{u} = (\Xb_{u,1} - \Xb_{u,2})/2$, and all the remaining terms are as
in (\ref{Eq:Model2}). The aim is thus at obtaining a design $\mathbf{u}\in\{
0, 1, 2\}^{N}$ so that $\Phi\{\Mb_u\}$ is minimized, where $\Mb_u =
\Fb
^T_{u}(\Ib_N - \omega\{[\mathbf{j}_N, \Eb_{u}]\})\Fb
_{u}$ and $\omega\{\Ab\}$
is the orthogonal projection matrix onto the space spanned by the
columns of the matrix $\Ab$. The following lemma can be easily proved.

%le3.1 #&#
\begin{lemma} \label{lemma6}
For a given design $\mathbf{u}\in\{0, 1, 2\}^N$, let $\bar
{\db}_u=(\bar
{d}_{u,1}, \ldots, \bar{d}_{u,N})^T \in\bar{\mathcal{D}} = \{-1,
0, 1\}
^N$ be defined as $\bar{d}_{u,n} = 0$, $1$ and $-1$ when $u_n = 0$, $1$
and $2$, respectively. Then
\[
\Mb_u \preceq\Fb^T_{u}\bigl(
\Ib_N - N^{-1}\Jb_N\bigr)\Fb_{u} =
\Xb^T_{\bar
{d}_u}\bigl(\Ib_N - N^{-1}
\Jb_N\bigr)\Xb_{\bar{d}_u}/4,
\]
where $\Xb_{\bar{d}_u} = [\bar{\db}_u, \Ub\bar{\db}_u, \ldots,
\Ub
^{K-1}\bar{\db}_u]$, and $\Ub$ is as in (\ref{Eq:U}).
In addition, $\Mb_u = \Xb^T_{\bar{d}_u}(\Ib_N - N^{-1}\Jb_N)\Xb
_{\bar
{d}_u}/4$ if $\mathbf{u}$ contains no zero.
\end{lemma}

Our approach for obtaining optimal fMRI designs for comparing HRFs is
by working on the upper bound of
$\Mb_u$ provided in Lemma~\ref{lemma6}. Specifically, we would like to
obtain a $\bar{\db}_u \in\bar{\mathcal{D}}$, or equivalently
a circulant CBWD $\Xb_{\bar{d}_u} \in\{-1, 0, 1\}^{N\times K}$, that
minimizes
$\Phi\{\Mb_b(\Xb_{\bar{d}_u})\}$.
As pointed out at the end of Section~\ref{sec2.2}, unlike the case
of estimating an HRF, here we also need to consider circulant CBWDs
with zero entries. If the obtained $\bar{\db}_u$
contains no zero, then the corresponding $\mathbf{u}$ is
$\Phi$-optimal. To
identify such a $\bar{\db}_u$, we consider the following lemma. For
convenience,
we omit the subscript of $\bar{\db}_u$ hereinafter, but its dependence
on $\mathbf{u}$ should be clear.

%le3.2 #&#
\begin{lemma} \label{lemma7}
Suppose $\bar{\db} \in\bar{\mathcal{D}}$ contains $r$ zeros, and
$\mathbf{j}
^T_N\bar{\db} = a$. We have the following results:
\begin{longlist}[(ii)]
\item[(i)] If $N = 4t - 1$, and $\Mb_b(\Xb_{\bar{d}}) = (N+1)[\Ib
_K -
N^{-1}\Jb_K]$, then
$a^2 = 1$, $r = 0$, and $\Mb(\Xb_{\bar{d}}) = \Xb^T_{\bar{d}}\Xb
_{\bar
{d}} = (N+1)\Ib_K - \Jb_K$.
\item[(ii)] If $N = 4t + 1$, and $\Mb_b(\Xb_{\bar{d}}) = (N-1)[\Ib
_K +
N^{-1}\Jb_K]$, then
$a^2 = 1$, $r = 0$, and $\Mb(\Xb_{\bar{d}}) = (N-1)\Ib_K + \Jb_K$.
\end{longlist}
\end{lemma}

\begin{pf}
We work only on (i) here. A similar argument can be applied to prove~(ii). For (i),
we have $\Mb_b(\Xb_{\bar{d}}) = \Mb(\Xb_{\bar{d}}) - (a^2/N) \Jb
_K =
(N+1)[\Ib_K - N^{-1}\Jb_K]$. Since each diagonal element of $\Mb(\Xb
_{\bar{d}})$ is an integer that is not greater than $N$, it can be seen
that $a^2 \leq1$. If $a = 0$, then $\Xb^T_{\bar{d}}\Xb_{\bar{d}} =
(N+1)[\Ib_K - N^{-1}\Jb_K]$. This leads to a contradiction since the
diagonal elements of the latter matrix are $(N+1)(1- N^{-1})$.
Therefore, $a^2 = 1$, $\Mb(\Xb_{\bar{d}}) = (N+1)\Ib_K - \Jb_K$
and $r
= 0$.
\end{pf}

The first main result in this section is an extension of
Theorem~\ref{thmm1}. We note that $\bar{N}_0(K, p_0)$ in Theorem~\ref
{thmm5} may not be
the same as $N_0(K, {p_0})$ in Theorem~\ref{thmm1}.

%th3.1 #&#
\begin{theorem} \label{thmm5}
Suppose $\bar{\db}^* \in\bar{\mathcal{D}}$ is a vector with $N = 4t-1$
elements, and it satisfies
$\Mb_b(\Xb_{\bar{d}^*}) = (N+1)[\Ib_K - N^{-1}\Jb_K]$. For any positive
number ${p_0}$,
there exists an $\bar{N}_0(K, {p_0})$ such that, if $N \geq\bar
{N}_0(K, {p_0})$,
then $\bar{\db}^*$ is $\Phi_p$-optimal over $\bar{\mathcal{D}}$
for any
$p \in[0, {p_0}]$.
\end{theorem}

\begin{pf}
Let $r$ be the number of zeros in $\bar{\db} \in\bar{\mathcal
{D}}$. It
is clear that
$\operatorname{tr}\{\Mb(\Xb_{\bar{d}})-N^{-1}\Jb_K\} = (N - r) - K/N$ is maximized
when $r = 0$. With Lemma~\ref{lemma7},
we can easily see that $\Mb(\Xb_{\bar{d}^*})-N^{-1}\Jb_K$ is
$(M,S)$-optimal over $\bar{\mathcal{M}}_b = \{\Mb(\Xb_{\bar{d}}) -
N^{-1}\Jb_K\}$. Following Lemma~\ref{Lemma_C1992} and Corollary~3.3 of
\citet{Cheng1987ar2021}, $\bar{\db}^*$ is $\Phi_p$-optimal over
$\bar{\mathcal{D}}$ for $p \in[0, {p_0}]$ when $N \geq\bar{N}_0(K,
{p_0})$ for some $\bar{N}_0(K, {p_0})$.
\end{pf}

An explicit lower bound $N_0(K, 1)$ for the $A$-criterion was given in
Theorem~\ref{4546465545485}. We show in the following theorem that one can take
$\bar{N}_0(K, 1)=N_0(K, 1)$.

%th3.2 #&#
\begin{theorem}\label{thmm5.1}
Let $\bar{\db}^* \in\bar{\mathcal{D}}$, where $N = 4t-1$, be such that
$\Mb_b(\Xb_{\bar{d}^*}) = (N+1)[\Ib_K - N^{-1}\Jb_K]$. If $K\geq
4$ and
$N \geq N_0(K, 1)$, where $N_0(K, 1)$ is as given in Theorem~\ref{4546465545485},
then $\bar{\db}^*$ is $A$-optimal (and $\Phi_p$-optimal for all $p
\in
[0,1]$) over $\bar{\mathcal{D}}$.
\end{theorem}

\begin{pf}
By Theorem~\ref{4546465545485}, it is enough to show that if $K\geq4$ and $N
\geq N_0(K, 1)$, then for any $\bar{\db}\in\bar{\mathcal{D}}$ that has
at least one zero entry, $\Phi_1\{\Mb_b(\Xb_{\bar{d}^*})\}\leq\Phi
_1\{
\Mb_b(\Xb_{\bar{d}})\}$. If $\bar{\db}\in\bar{\mathcal{D}}$ has at
least one zero entry, then each diagonal entry of $\Mb_b(\Xb_{\bar
{d}})$ is at most $N-1$; thus $\Phi_1\{\Mb_b(\Xb_{\bar{d}})\}\geq
K/(N-1)$. On the other hand, since $\Mb_b(\Xb_{\bar{d}^*})=(N+1)[\Ib_K
- N^{-1}\Jb_K]$ has two distinct eigenvalues $N+1$ and $(N+1)(N-K)/N$,
with multiplicity $K-1$ and 1, respectively, we have $\Phi_1\{\Mb
_b(\Xb
_{\bar{d}^*})\}=(K-1)/(N+1)+N/[(N+1)(N-K)]$. It follows that $\Phi_1\{
\Mb_b(\Xb_{\bar{d}^*})\}\leq\Phi_1\{\Mb_b(\Xb_{\bar{d}})\}$ provided
$(K-1)/(N+1)+N/[(N+1)(N-K)]\leq K/(N-1)$. The latter is the same as
$N\geq2K-1$. Therefore, it remains to show that $2K-1\leq N_0(K, 1)$.
Since $N_0(K, 1)$ is the greatest real root of the cubic function
$c(x) = 2x^3 + (10-7K)x^2 + 2(2K-5)(K-1)x + 4K^2 - 7K$, $c(x)>0$ for
all $x>N_0(K, 1)$. One can verify that $c(2K-1)<0$ when $K\geq4$. It
follows that in this case $2K-1<N_0(K, 1)$.
\end{pf}

For $N = 4t$, circulant OAs or equivalently circulant partial Hadamard
matrices described in Section~\ref{sec2.3} can be used to construct
$\bar{\db}^*$
that has $\Mb_b(\Xb_{\bar{d}^*}) = N\Ib_{K}$. Such a $\bar{\db
}^*$ can
be easily seen to be universally optimal in $\bar{\mathcal{D}}$.

%Our
Theorem~\ref{thmm6} below
%also
helps to identify some optimal $\bar{\db}$ for $N = 4t + 1$.
For deriving this theorem, we again consider Lemmas \ref{C1980} and
\ref{lemma5}.

%th3.3 #&#
\begin{theorem}\label{thmm6}
For $N = 4t+1$, let $\bar{\db}^* \in\bar{\mathcal{D}}$ have
$\Mb_b(\Xb_{\bar{d}^*}) = (N-1)[\Ib_K + N^{-1}\Jb_K]$.
Then, $\bar{\db}^*$ is optimal over $\bar{\mathcal{D}}$ for any
type 1
criterion.
\end{theorem}

\begin{pf}
$\Mb_b(\Xb_{\bar{d}^*})$ has two nonzero eigenvalues, and the smaller
eigenvalue has multiplicity $K-1$. It can also be seen that condition
(a) in Lemma~\ref{C1980} is satisfied by $\bar{\db}^*$. In addition,
$\operatorname{tr}\{\Mb_b(\Xb_{\bar{d}^*})\} = K(N - N^{-1})$, and $\operatorname{tr}\{\Mb
^2_b(\Xb
_{\bar{d}^*})\} - (\operatorname{tr}\{\Mb_b(\Xb_{\bar{d}^*})\})^2/K
= K(K-1)(1 - N^{-1})^2$. Thus, condition (b) of Lemma~\ref{C1980}
is satisfied if and only if
%
%e3.3 #&#
\begin{equation}\qquad
K\bigl(N - N^{-1}\bigr) - A_{\bar{d}} \geq\sqrt{
\frac{K}{K-1}} \biggl[\bigl(1 - N^{-1}\bigr)\sqrt{K(K-1)} - \sqrt
{B_{\bar{d}} - \frac{A^2_{\bar{d}}}{K}} \biggr] \label{Eq:condc}
\end{equation}
for $\bar{\db} \in\bar{\mathcal{D}}$, where $A_{\bar{d}} = \operatorname{tr}\{
\Mb
_b(\Xb_{\bar{d}})\}$, and
$B_{\bar{d}} = \operatorname{tr}\{\Mb^2_b(\Xb_{\bar{d}})\}$. Clearly, (\ref{Eq:condc})
holds for the class of $\bar{\db} \in\bar{\mathcal{D}}$ that satisfy
$A_{\bar{d}} \leq K(N - N^{-1}) -\sqrt{K/(K-1)}(1 - N^{-1})\sqrt
{K(K-1)} = K(N-1)$. Thus, all $\bar{\db}$'s in this class are
outperformed by $\bar{\db}^*$ with respect to any type 1 criterion. For
any other $\bar{\db}$, let $r$ be the number of zeros
in $\bar{\db}$ and $a=\mathbf{j}_N^T \bar{\db}$. Then we
have $A_{\bar{d}} >
K(N-1)$, and $(1 - r - a^2/N) > 0$ since $A_{\bar{d}} = K[(N-r) -
a^2/N]$. Consequently, $r = 0$, and $\bar{\db}$ contains no zero.
Following Lemma~\ref{lemma5}, $\bar{\db}^*$ is also optimal over the
class of designs with no zero for any type 1 criterion. Our claim then follows.
\end{pf}

With these results, we can derive the following theorem for identifying
some optimal fMRI designs for studying contrasts between two HRFs.

%th3.4 #&#
\begin{theorem}\label{thmm7}
Suppose $\db_H$ is a Hadamard sequence, $\db_{1,g, H}$ is defined as in
Theorem~\ref{thmm3}, $\db_{2,g, H}$ is
as in Theorem~\ref{thmm4}, and $\Mb_u$ is the information matrix for
$\zetab$ in model (\ref{Eq:Model3}) for a design $\mathbf{u}\in\{0, 1, 2\}
^N$. We have the following results:
\begin{longlist}[(a)]
\item[(a)] Suppose $N = 4t - 1$, and $\mathbf{u}^* = \mathbf{j}_N + \db_H$ or $\mathbf{u}^*
= 2\mathbf{j}_N - \db_H$. If ${p_0} >0$, $p \in[0,
{p_0}]$, and
$N \geq\bar{N}_0(K, {p_0})$ for some $\bar{N}_0(K, {p_0}) > 0$, then
$\Phi_p\{\Mb_{u^*}\} \leq\Phi_p\{\Mb_u\}$ for any $\mathbf{u}\in\{0, 1, 2\}
^N$; here, $\Phi_p$ is defined as in Definition~\ref{Def:Phip}. For
$K\geq4$, we can take $\bar{N}_0(K, 1) $ to be the $N_0(K, 1)$ given
in Theorem~\ref{4546465545485}.
\item[(b)] Suppose $N = 4t$, and $\mathbf{u}^* = \mathbf{j}_N + \db_{1,g,H}$ or $\mathbf{u}
^* = 2\mathbf{j}_N - \db_{1,g,H}$. If $K \leq g + 1$, and
$\Phi$ is any
criterion satisfying the conditions in Definition~\ref{Def:UO}, then
$\Phi\{\Mb_{u^*}\} \leq\Phi\{\Mb_u\}$ for any $\mathbf{u}\in\{0, 1, 2\}^N$.
\item[(c)] Suppose $N = 4t +1$, and $\mathbf{u}^* = \mathbf{j}_N + \db_{2,g,H}$ or
$\mathbf{u}^* = 2\mathbf{j}_N - \db_{2,g,H}$. If
$K \leq g + 1$, and $\Phi_{(f)}$ is
any type 1 criterion defined in Definition~\ref{Def:Type1}, then $\Phi
_{(f)}\{\Mb_{u^*}\} \leq\Phi_{(f)}\{\Mb_u\}$ for any $\mathbf{u}\in\{0, 1,
2\}^N$.
\end{longlist}
\end{theorem}

\begin{pf}
For all the designs $\mathbf{u}^*$ in (a), (b) and (c), we
have $\Mb_{u^*} =
\Xb^T_{\bar{d}_{u^*}}(\Ib_N - N^{-1}\Jb_N)\Xb_{\bar{d}_{u^*}}/4$, where
$\bar{\db}_{u^*}$ is defined as in Lemma~\ref{lemma6}. When $\mathbf{u}^* = \mathbf{j}
_N + \db_H$ (or, resp., $\mathbf{u}^* = 2\mathbf{j}_N - \db_H$), we have $\Xb
_{\bar{d}_{u^*}} = \Xb_{\bar{d}_H}$ with $\bar{\db}_H = \mathbf{j}_N - 2\db_H$
(or, resp., $\bar{\db}_H = 2\db_H - \mathbf{j}_N$). We thus have that,
if $N \geq\bar{N}_0(K, {p_0})$ and $p \in[0, {p_0}]$, then
\[
\Phi_p\{\Mb_{u^*}\} = \Phi_p\bigl\{
\Mb_b(\Xb_{\bar{d}_H})/4\bigr\} \leq \Phi_p\bigl\{
\Mb_b(\Xb_{\bar{d}_u})/4\bigr\} \leq\Phi_p\{
\Mb_u\}
\]
for any $\mathbf{u}\in\{0, 1, 2\}^N$. This completes the
proof for (a).
Similar arguments can be used to prove (b) and (c) and are omitted.
\end{pf}

%s4 #&#
\section{Conclusion}\label{sec4}
Neuroimaging experiments utilizing the pioneering fMRI technology are
widely conducted in a variety of research fields for gaining better knowledge
about human brain functions. One of the key steps to ensure the success
of such an experiment is to judiciously select an optimal fMRI design.
Existing studies on obtaining optimal fMRI designs primarily focus on
computational approaches. However, insightful analytical results, while
important, are
rather scarce and scattered. To address this important issue, we
conduct a systematic and analytical analysis to characterize some
optimal fMRI designs for
estimating the HRF of a stimulus type and for comparing HRFs of two
stimulus types. Under certain conditions, we show that the popularly
used binary $m$-sequences
as well as the more general Hadamard sequences are optimal in some
statistically meaningful senses. We also identify several new classes
of high-quality fMRI designs
and present systematic methods for constructing them. These designs
exist in many design lengths where good fMRI designs have not been
reported previously.
There, however, are many research challenges that need to be overcome.
For example, our results provide good designs for design lengths of
$N = 4t - 1$, $4t$ and $4t+1$. A~future research of interest is on
identifying optimal fMRI designs for cases with $N = 4t + 2$. In
addition, our experience
indicates that the designs that we present here remain quite efficient
under some violations of model assumptions [cf. \citet{Kao2014ar51}].
Nevertheless, it still is of interest to analytically study optimal
designs for other situations (e.g., with an autoregressive error
process). Extending current results to
cases with a greater number of stimulus types is also a future research
of interest. Many research opportunities exist in this new and
wide-open research area.

\begin{appendix}\label{app}
% redefine the command that creates the equation no.
%\setcounter{equation}{0} % reset counter
%\setcounter{definition}{0} % reset counter
%\setcounter{lemma}{0} % reset counter

\section*{Appendix: A proof of Theorem~\texorpdfstring{\protect\ref{4546465545485}}{2.2}}
For $N = 3\ (\mathrm{mod}\ 4) \geq4$ and $K \geq4$, we consider the following
set of $K$-by-$K$ nonnegative definite matrices:
\[
\Xi_{N,K} = \bigl\{ \Eb_K = \bigl((e_{ij})
\bigr)_{i,j = 1, \ldots, K} \mid e_{ij} = 3\ (\mathrm{mod}\ 4)\ \forall i,j,
e_{ii} = N, \Eb_K \succ N^{-1}\Jb_K
\bigr\},
\]
where $\Eb_K \succ N^{-1}\Jb_K$ indicates that $\Eb_K - N^{-1}\Jb
_K$ is
positive definite. With Lemma~\ref{Lemma_mod3}, it can be seen that
$\Mb(\Xb_{\tilde{d}}) \in\Xi_{N,K}$ for any $\tilde{\db} \in
\tilde
{\mathcal{D}}$ having a nonsingular
$\Mb_b(\Xb_{\tilde{d}})$. The idea for proving Theorem~\ref{4546465545485} is
then to show that
an $\Eb_K \in\Xi_{N,K}$ minimizing $\operatorname{tr}\{[\Eb_K - N^{-1}\Jb
_K]^{-1}\}$
is similar (with some permutations of rows and columns)
to a block matrix $\Bb\in\Xi_{N,K}$ to be defined in Definition~\ref
{DefA1} below. We also will show in Lemma~\ref{lemmaA3} that, when the
condition in Theorem~\ref{4546465545485} is satisfied,
we have $\operatorname{tr}\{[\Mb(\Xb_{\tilde{d}^*}) - N^{-1}\Jb_K]^{-1}\} = \min_{B
\in\mathcal{B}} \operatorname{tr}\{[\Bb- N^{-1}\Jb_K]^{-1}\}$, where
$\mathcal{B} \subset\Xi_{N,K}$ is the set of all block matrices. With
these facts and Lemma~\ref{Lemma_mod3}, we have
%
%e4.1 #&#
\begin{eqnarray}\label{Eq:phi1}
&&\Phi_1\bigl\{\Mb_b(\Xb_{\tilde{d}^*})\bigr\}\nonumber\\
&&\qquad=
\Phi_1\bigl\{\Mb(\Xb_{\tilde
{d}^*}) - N^{-1}
\Jb_K \bigr\} = \min_{E_k \in\Xi_{N,K}} \Phi_1\bigl
\{\Eb_K - N^{-1}\Jb_K \bigr\}
\\
&&\qquad\leq\min_{\tilde{d} \in\tilde{\mathcal{D}}} \Phi_1\bigl\{\Mb(\Xb
_{\tilde
{d}}) - N^{-1}\Jb_K \bigr\} \leq \min
_{\tilde{d} \in\tilde{\mathcal{D}}} \Phi_1\bigl\{\Mb_b(\Xb
_{\tilde
{d}})\bigr\}. \nonumber
\end{eqnarray}
Our claim in Theorem~\ref{4546465545485} then follows from (\ref{Eq:phi1}), and
Corollary~3.3 of \citet{Cheng1987ar2021}. This approach is similar
to that of \citet{SatheShenoy1989ar2021}, and \citet
{GalilKiefer1980ar2021}, where
weighing designs under the unbiasedness assumption are considered. We
now present the details
of our proof.

%de4.1 #&#
\begin{definition}\label{DefA1}
A block matrix $\Bb\in\Xi_{N,K}$ is of the form
\[
\Bb= \bigoplus_{i = 1}^m \bigl[(N-3)
\Ib_{r_i} + 4 \Jb_{r_i} \bigr] - \Jb_K
\]
for an integer $m \in\{1, \ldots, K\}$ representing the number of
``blocks.'' Here, $\bigoplus$ is the matrix direct sum, and $r_1, r_2,
\ldots, r_m$ are the block sizes of $\Bb$ that satisfy $r_i \geq1$, and
$\sum_{i=1}^m r_i = K$.
\end{definition}

For these block matrices, the following result is an extension of
Theorem~2.1(a) of \citet{Masaro1988ar2021} and
equation (1.1) of \citet{SatheShenoy1989ar2021}. Clearly, this
result also implies that
$\operatorname{tr}\{\Bb^{-1}_b\}$ is invariant to a rearrangement of the block sizes
$r_1, \ldots, r_m$, which
facilitates the derivation of the subsequent results; here $\Bb_b =
\Bb
- N^{-1}\Jb_K$.
% for $\Bb\in\mathcal{B}$.

%le4.1 #&#
\begin{lemma} \label{lemmaA1}
Let $\mathcal{B} \subset\Xi_{N,K}$ be the set of all block matrices. Then
%, and $\Bb_b = \Bb- N^{-1}\Jb_K$ f
for $\Bb\in\mathcal{B}$,
we have
%
%e4.2 #&#
\begin{eqnarray}\label{TwoEqu}
\operatorname{tr}\bigl\{\Bb^{-1}_b\bigr\} &=& \sum
_{i=1}^m L_i^{-1} +
\frac{K-m}{N - 3} + \frac
{\sum_{i=1}^m r_i L_i^{-2} }{(1+N^{-1})^{-1} - \sum_{i=1}^m r_i
L_i^{-1}}
\nonumber
\\[-8pt]
\\[-8pt]
\nonumber
&=&\frac{1}{4t} \Biggl\{K - \sum_{i=1}^m
\frac{r_i}{t + r_i} + \frac{t
\sum_{i=1}^m {r_i}/{(t + r_i)^2}}{{4}/{(1+N^{-1})}
- \sum_{i=1}^m {r_i}/{(t + r_i)}} \Biggr\},
\end{eqnarray}
where $L_i = N - 3 + 4 r_i$, and $t = (N-3)/4 \geq1$.
\end{lemma}

\begin{pf} Let $\Cb_{r_i} = (N-3)\Ib_{r_i} + 4 \Jb_{r_i}$, we have
$\Bb_b = \bigoplus_{i = 1}^m \Cb_{r_i} - (1 + N^{-1})\Jb_K$, and
\begin{eqnarray*}
\Bb^{-1}_b &=& \Biggl[ \bigoplus
_{i = 1}^m \Cb_{r_i} - \bigl(1 +
N^{-1}\bigr)\mathbf{j} _K\mathbf{j}'_K \Biggr]^{-1}
\\
&=&\bigoplus_{i = 1}^m
\Cb^{-1}_{r_i} + \frac{(1 + N^{-1})
[\bigoplus_{i = 1}^m \Cb^{-1}_{r_i} ]\mathbf{j}_K\mathbf{j}'_K [\bigoplus_{i
= 1}^m \Cb^{-1}_{r_i} ]}{1 - (1 + N^{-1}) \sum_{i=1}^m \mathbf{j}_{r_i}
\Cb^{-1}_{r_i}\mathbf{j}_{r_i}}.
\end{eqnarray*}
The two equalities in (\ref{TwoEqu}) can then be derived by some
simple algebra.
\end{pf}

The following lemma indicates that a block matrix minimizing the trace
of $\Bb_b^{-1} = [\Bb- N^{-1}\Jb_K]^{-1}$ can be found over a small
subset of $\mathcal{B}$
[cf. Theorem~2.1(b) of \citet{Masaro1988ar2021}].

%le4.2 #&#
\begin{lemma} \label{lemmaA2}
Let $\mathcal{B}_s \subset\mathcal{B}$ be the set of block matrices
having blocks of only one size
or two contiguous sizes. Then
\[
\min_{\Bb\in\mathcal{B}_s} \operatorname{tr}\bigl\{ \Bb_b^{-1} \bigr
\} = \min_{
\Bb
\in\mathcal{B}} \operatorname{tr}\bigl\{ \Bb^{-1}_b
\bigr\}.
\]
\end{lemma}

\begin{pf}
Among the block matrices that yield the minimum $\operatorname{tr}\{ \Bb^{-1}_b \}$
over $\mathcal{B}$, let $\Bb_m$ be the one with the smallest number of
blocks, $m$. Clearly, we only need to consider cases where $m \geq2$.
Without loss of generality (see also the statement above Lemma~\ref
{lemmaA1}), we may assume that the first two block sizes $r_1$ and
$r_2$ are, respectively, the largest and the smallest block sizes among
the $m$ block sizes.
With Lemma~\ref{lemmaA1}, we can then write
$\operatorname{tr}\{\Bb_{m,b}^{-1}\} = [K + f(x)]/4t$, where $\Bb_{m,b} = \Bb_m -
N^{-1}\Jb_K$,
\begin{eqnarray*}
f(x) &=& \frac{t [{x}/{(t+x)^{2}}+{(r-x)}/{(t+r-x)^2} +
\beta
 ]}{
{4}/{(1+N^{-1})} -  [ {x}/{(t + x)} + {(r-x)}/{(t+r-x)} +
\alpha ]} \\
&&{}- \biggl[\frac{x}{t + x} + \frac{r-x}{t+r-x} +
\alpha \biggr],
\end{eqnarray*}
$r = r_1 + r_2$, $x = r_1$ (or $r_2$) with $0 < x < r$, $\alpha= \sum_{i = 3}^{m} r_i/(t+r_i)$, and
$\beta= \sum_{i = 3}^{m} r_i/(t+r_i)^2$. We note that this expression
of $\operatorname{tr}\{\Bb_b^{-1}\}$
applies to any block matrix, and that $\alpha= \beta= 0$ for block
matrices with two or fewer blocks. Since the
number of blocks of $\Bb_m$ is $m \geq2$, we have $\operatorname{tr}\{ \Bb^{-1}_{m,b}
\} < \operatorname{tr}\{ \Bb^{-1}_b \}$ for those block matrices $\Bb$ with only one
block; thus,
$f(x) < f(0) = f(r)$. With some simple algebra similar to that of
\citet
{Masaro1988ar2021}, we also have
\[
f(x) = h(y) = \frac{ay^2 + b}{c y^2 + d}\quad \mbox{and}\quad f'(x) =
h'(y) = \frac{2(ad-bc)y}{(cy^2 + d)^2},
\]
where $y = x - 0.5 r \in(-0.5r, 0.5r)$, and $a, b, c, d$ are some
constants. Along with the fact that $f(x) < f(0) = f(r)$, and $f$ is
symmetric about $0.5r$, the minimum of $f(x)$ occurs when $x$ is the
integer closest to $0.5r$. Consequently, $r_1 = r_2$ when $r$ is even, and
$r_1 = r_2 +1$ when $r$ is odd.
\end{pf}

With these results, we can now work on $\mathcal{B}_{0,s} \subset
\mathcal{B}_s$ that consists of block matrices having
$m (< K)$ block sizes with $r_1 -1 = \cdots= r_v - 1 = r_{v+1} =
\cdots= r_{m} = r \geq1$, and $v \geq1$. We note that,
for any $m_0 \geq1$ and $r_0 \geq2$, a block matrix having $(m, v, r)
= (m_0, 0, r_0)$ can be treated as a block matrix with
$(m, v, r) = (m_0, m_0, r_0 -1)$, and is thus in $\mathcal{B}_{0,s}$;
see also \citet{SatheShenoy1989ar2021}. Consequently,
the only block matrix in $\mathcal{B}_s$ that is left out from
$\mathcal
{B}_{0,s}$ is $\Bb^* = (N+1)\Ib_K - \Jb_K$, which has
$(m, v, r) = (K, 0, 1)$, or equivalently, $(m, v, r) = (K, K, 0)$.
Under the condition described in the following lemma, we have $\operatorname{tr}\{
(\Bb
^* - N^{-1}\Jb_K)^{-1} \} \leq \operatorname{tr}\{ (\Bb_s - N^{-1}\Jb_K)^{-1} \}$
for any other $\Bb_s \in\mathcal{B}_s$.

%le4.3 #&#
\begin{lemma} \label{lemmaA3}
Let $\Bb^* = (N+1)\Ib_K - \Jb_K$, $\Bb_s \in\mathcal{B}_{0,s}$ be
a
previously described block matrix,
$\Bb^*_{b} = \Bb^* - N^{-1}\Jb_K$, $\Bb_{s,b} = \Bb_{s} -
N^{-1}\Jb_K$,
and $N_0(K, 1)$ be the greatest real root of the cubic function
$c(x) = 2x^3 + (10-7K)x^2 + 2(2K-5)(K-1)x + 4K^2 - 7K$. If $N \geq
N_0(K, 1)$, then
\[
\operatorname{tr}\bigl\{ \bigl(\Bb_b^*\bigr)^{-1} \bigr\} \leq \operatorname{tr}\bigl\{
\Bb_{s,b}^{-1} \bigr\}.
\]
\end{lemma}

\begin{pf}
Let $u = m-v$, and $\Bb$ be obtained by replacing a block of size $r +
1$ in $\Bb_s$ with a block of size $r$
and a block of size $1$. From (\ref{TwoEqu}), we have
\[
\operatorname{tr}\bigl\{\Bb_{s,b}^{-1} \bigr\} = \frac{K-m}{N - 3} +
\frac{u - 1}{L} + \frac
{v -
1}{L + 4} + \frac{(2L + 4)(1+N^{-1})^{-1} - K}{g(v)},
\]
where $L = N - 3 + 4r$, and $g(v) = 4v(r+1) + (L+4)\xi$, and $\xi=
[L(1+N^{-1})^{-1}-K]$. In addition,
\begin{eqnarray*}
\operatorname{tr}\bigl\{\Bb_b^{-1} \bigr\} &=& \frac{K-m-1}{N - 3} +
\frac{u}{L} + \frac{v -
2}{L + 4} + \frac{1}{N+1}
\\
&&{}+ \frac{(2L + 4)(1+N^{-1})^{-1} - K + 16r(r-1)(N+1)^{-2}}{g(v) -
4r(N+1+L)(N+1)^{-1}}.
\end{eqnarray*}
Since $g(v) \geq g(1)$, we have
\begin{eqnarray*}
&&\operatorname{tr}\bigl\{\Bb_{s,b}^{-1} - \Bb_b^{-1}
\bigr\}
\\
&&\qquad = \frac{4r}{L(N-3)} + \frac{(2L + 4)(1+N^{-1})^{-1} - K}{g(v)} +\frac
{1}{L+4} -
\frac{1}{N+1}
\\
&&\qquad\quad{} - \frac{(2L + 4)(1+N^{-1})^{-1} - K + 16r(r-1)(N+1)^{-2}}{g(v) -
4r(N+1+L)(N+1)^{-1}}
\\
&&\qquad\geq\frac{4r}{L(N-3)} + \frac{(2L + 4)(1 + N^{-1})^{-1} - K}{g(1)}
\\
&&\qquad\quad{} - \frac{(2L+4 - 4r)(1+N^{-1})^{-1} - K}{(N+1)\xi-4(r-1)} = \Delta (r). \label{Eq:Ineq}
\end{eqnarray*}
The equality holds when $v = 1$. With some algebra,
we have
\begin{eqnarray*}
\Delta(r) &=& \frac{8r(N+12r-11)(N+1)^{-1}c(N)}{L(N-3)[(N+1)\xi-
4(r-1)][(L+4)\xi+ 4(r+1)]}
\\
&&{} + 16r(r-1) \\
&&\hspace*{-2pt}\quad{}\times\frac{N^2 - 7N + 12 - 8(1+N^{-1})^{-1} +
16(r-1)^2N(3N-1)(N+1)^{-2}}{L(N-3)[(N+1)\xi- 4(r-1)][(L+4)\xi+
4(r+1)]}
\\
&&{} + 16r(r-1) \\
&&\hspace*{-2pt}\quad{}\times\frac{4(r-1)(N+1)^{-1}[(7N-K)(N-K) + 5(N^2 - 1) +
8N]}{L(N-3)[(N+1)\xi- 4(r-1)][(L+4)\xi+ 4(r+1)]},
\end{eqnarray*}
where $c(N) = 2N^3 + (10-7K)N^2 + 2(2K-5)(K-1)N + 4K^2 - 7K$.
With $N = 3\ (\mathrm{mod}\ 4) \geq4$ and $K \geq4$, it can be seen that, when
$N \geq N_0(K, 1)$, we have $c(N) \geq0$, $\Delta(r) > 0$ for $r > 1$,
and $\Delta(1) \geq0$. Consequently, for a $\Bb_s$, we either (i) find
a $\Bb\notin\mathcal{B}_s$ with
$\operatorname{tr}\{\Bb_{s,b}^{-1}\} > \operatorname{tr} \{\Bb_b^{-1}\}$, or (ii) can keep splitting
each block of size $r+1 = 2$ in $\Bb_s$
into two blocks of size $1$ without increasing the objective function
($\operatorname{tr}\{\Bb_b^{-1}\}$).
For the first case, $\Bb_s (\neq\Bb^*)$ is obviously not an optimal
block matrix
of our interest. For the second case, we can continue the process until
$\Bb= \Bb^*$. Our claim thus follows.
\end{pf}

The results we have so far suggest that, under the condition of Lemma~\ref{lemmaA3}, $\Bb^* = (N+1)\Ib_K - \Jb_K$ minimizes $\operatorname{tr}\{\Bb
_b^{-1}\}
$ over all block matrices $\Bb$. \red{With the following
Lemma~\ref{lemmaA4}, our proof of Theorem~\ref{4546465545485} is then complete.}

%le4.4 #&#
\begin{lemma} \label{lemmaA4}
Let $\Bb^* = (N+1)\Ib_K - \Jb_K$. If $N \geq N_0(K, 1)$ for the $N_0(K,
1)$ defined in Lemma~\ref{lemmaA3}, then
\[
\operatorname{tr}\bigl\{\bigl(\Bb_b^*\bigr)^{-1}\bigr\} = \min
_{\Eb_K \in\Xi_{N,K}} \operatorname{tr}\bigl\{\bigl[\Eb_K - N^{-1}\Jb
_K\bigr]^{-1}\bigr\}.
\]
\end{lemma}

The proof of Lemma~\ref{lemmaA4} is lengthy, but otherwise is a
simple extension of that of Theorem~2.2 of \citet{SatheShenoy1989ar2021}.
The main idea is to show that an $\Eb^*_K \in\Xi_{N,K}$
minimizing $\operatorname{tr}\{(\Eb_{b,K}^{-1}\}$ is similar to a block matrix after
some permutations of rows and columns. Lemma~\ref{lemmaA4}
then follows from Lemma~\ref{lemmaA3}. To that end, we need the
following lemmas, which are extensions of
results in \citet{SatheShenoy1989ar2021}. Lemma~\ref{lemmaS1} is
a well-known result, and
the proof is omitted. We also use the following notation: $\Eb^*_K \in
\Xi_{N,K}$ is a
matrix such that $\operatorname{tr}\{(\Eb^*_{b,K})^{-1}\} = \min_{\Eb_K \in\Xi_{N,K}}
\operatorname{tr}\{\Eb^{-1}_{b,K}\}$,
$N_b = N - N^{-1}$, $3_b = 3 - N^{-1}$, $c_b = c - N^{-1}$ for some $c
= 3\ (\mathrm{mod}\ 4)$, and $\mub_{b,i} = \mub_i - N^{-1} \mathbf{j}_{K-2}$ and $\mub
_{b,j}= \mub_j - N^{-1} \mathbf{j}_{K-2}$ for some $\mub
_i$ and $\mub_j$ whose
elements are congruent to 3 modulo 4. In addition, $a_{b,i,j} = \mub
^T_{b,i}\Eb^{-1}_{b,K-2}\mub^T_{b,j}$, $b_{b, i, j} = \mub
^T_{b,i}\Eb
^{-2}_{b,K-2}\mub^T_{b,j}$, $A_{b, i, j} = N_b - a_{b, i, j}$,
$z_{b,i,j}(c_b) = c_b - a_{b,i,j}$,
and
%
%e4.3 #&#
\begin{eqnarray}\label{Eq:fij}
&&f_{b,i,j}(c_b)
\nonumber
\\[-8pt]
\\[-8pt]
\nonumber
&&\qquad = \frac{(A_{b,i,i} + A_{b,j,j}) + A_{b,i,i}
b_{b,j,j}+A_{b,j,j}b_{b,i,i} - 2b_{b,
i,j}z_{b,i,j}(c_b)}{A_{b,i,i}A_{b,j,j} - z^2_{b,i,j}(c_b)}.
\end{eqnarray}

%le4.5 #&#
\begin{lemma}\label{lemmaS1}
Let $\Eb= ((\Eb_{ij}))$ for $i, j = 1, 2$ be a partitioned positive
definite matrix, where
$\Eb_{11}$ and $\Eb_{22}$ are square matrices. We have
\[
\operatorname{tr}\bigl\{\Eb^{-1}\bigr\} = \operatorname{tr}\bigl\{\Eb^{-1}_{22}
\bigr\} + \operatorname{tr}\bigl\{\Vb\bigl[\Ib+ \Eb_{12}\Eb _{22}^{-2}
\Eb_{21}\bigr]\bigr\},
\]
where $\Vb= (\Eb_{11} - \Eb_{12}\Eb_{22}^{-1}\Eb_{21})^{-1}$.
\end{lemma}

%le4.6 #&#
\begin{lemma}\label{lemmaS2}
$\operatorname{tr}\{(\Eb^*_{b,K})^{-1}\} <
 \operatorname{tr}\{(\Eb^*_{b,K-1})^{-1}\} + (N-3)^{-1}$.
\end{lemma}

\begin{pf}
For $K = 2$, $\operatorname{tr}\{\Eb^{-1}_{b, 2}\} = 2N_b/(N_b^2 - c_b^2)$ is
minimized when $c = -1$, or equivalently, $c_b = -1 - N^{-1}$.
Thus, $\operatorname{tr}\{(\Eb^*_{b, 2})^{-1}\} - \operatorname{tr}\{(\Eb^*_{b, 1})^{-1}\} =$
\[
\frac{2N_b}{N_b^2 - (1+N^{-1})^2} - N_b = \frac
{N^2-2N+2}{(N+1)(N-1)(N-2)} < \frac{1}{N-3}.
\]
Suppose $\operatorname{tr}\{(\Eb^*_{b,K-1})^{-1}\} < \operatorname{tr}\{(\Eb^*_{b,K-2})^{-1}\} +
(N-3)^{-1}$.
We would like to show that $\operatorname{tr}\{(\Eb^*_{b,K})^{-1}\} < \operatorname{tr}\{(\Eb
^*_{b,K-1})^{-1}\} + (N-3)^{-1}$. To that end, we write
\[
\Eb^*_{b,K-1} = \left[ %
\matrix{
N_b & \mub^T_{b,j}
\vspace*{2pt}\cr
\mub_{b,j} & \Eb_{b,K-2} }
 \right].
\]
With Lemma~\ref{lemmaS1} and the fact that $\Eb^*_{b,K-1}$ is positive
definite, we have
$b_{b,j,j} > 0$ and
\begin{eqnarray*}
A^{-1}_{b,j,j} &<& A^{-1}_{b,j,j}(1+b_{b,j,j})
= \operatorname{tr}\bigl\{\bigl(\Eb^*_{b,K-1}\bigr)^{-1}\bigr\} - \operatorname{tr}\bigl\{
\Eb_{b,K-2}^{-1}\bigr\}
\\
& \leq &\operatorname{tr}\bigl\{\bigl(\Eb^*_{b,K-1}\bigr)^{-1}\bigr\} - \operatorname{tr}
\bigl\{\bigl(\Eb^*_{b,K-2}\bigr)^{-1}\bigr\} <
(N-3)^{-1}.
\end{eqnarray*}
Thus, $A_{b, j, j} = N_b - a_{b, j, j} > N - 3$, and this in turn implies
$z_{b,j,j}(3_b) = 3_b - a_{b, j, j} > 0$. Following this fact and some
simple algebra, we can show that
the following matrix $\Eb_{b,K}(j)$, which is obtained by adding a row
and a column to
$\Eb^*_{b,K-1}$, is positive definite, and is thus in $\Xi_{N,K}$:
\[
\Eb_{b,K}(j) = \left[ %
\matrix{
N_b & 3_b & \mub^T_{b,j}
\vspace*{2pt}\cr
3_b & N_b & \mub^T_{b,j}
\vspace*{2pt}\cr
\mub_{b,j} & \mub_{b,j} & \Eb_{b,K-2} }
 \right].
\]
With Lemma~\ref{lemmaS1}, we also have
%
%e4.4 #&#
\begin{eqnarray}\label{Eq:Ej}
\operatorname{tr}\bigl\{\Eb^{-1}_{b,K}(j)\bigr\} &=& \operatorname{tr}\bigl\{
\Eb_{b,K-2}^{-1}\bigr\} + f_{b,j,j}(3_b)
\nonumber
\\[-8pt]
\\[-8pt]
\nonumber
&=&\operatorname{tr}\bigl\{\bigl(\Eb^*_{b,K-1}\bigr)^{-1}\bigr\} -
A^{-1}_{b,j,j}(1+b_{b,j,j}) + f_{b,j,j}(3_b).
\nonumber
\end{eqnarray}
By noting that $f_{b,j,j}(3_b)$ in (\ref{Eq:fij}) can be written as
\[
f_{b,j,j}(3_b) = \frac{2[z_{b,j,j}(3_b) +
(N-3)(1+b_{b,j,j})]}{(N-3)(N-3 + 2z_{b,j,j}(3_b))},
\]
we can show that $f_{b,j,j}(3_b) - A^{-1}_{b,j,j}(1+b_{b,j,j}) < (N-3)^{-1}$.
The proof is then completed by the fact that
\begin{eqnarray*}
\operatorname{tr}\bigl\{\bigl(\Eb^*_{b,K}\bigr)^{-1}\bigr\} &\leq& \operatorname{tr}\bigl\{
\Eb^{-1}_{b,K}(j)\bigr\} =\operatorname{tr}\bigl\{\bigl(\Eb
^*_{b,K-1}\bigr)^{-1}\bigr\} + f_{b,j,j}(3_b)
- A^{-1}_{b,j,j}(1+b_{b,j,j})
\\
&<& \operatorname{tr}\bigl\{\bigl(\Eb^*_{b,K-1}\bigr)^{-1}\bigr\} +
(N-3)^{-1}.
\end{eqnarray*}
\upqed\end{pf}

%le4.7 #&#
\begin{lemma}\label{lemmaS3}
Write $\Eb^*_{b,K} = \Eb^*_K - N^{-1}\Jb_K$ in the form of
%
%e4.5 #&#
\begin{equation}
\Eb^*_{b,K} = \left[\matrix{
N_b & c_b & \mub^T_{b,i}
\vspace*{2pt}\cr
c_b & N_b & \mub^T_{b,j}
\vspace*{2pt}\cr
\mub_{b,i} & \mub_{b,j} & \Eb_{b,K-2}}
 \right] = \left[ \matrix{ N & c &
\mub^T_{i}
\vspace*{2pt}\cr
c & N & \mub^T_{j}
\vspace*{2pt}\cr
\mub_{i} & \mub_{j} & \Eb_{K-2}}
 \right] - N^{-1}\Jb_K. \label{Eq:Estar}
\end{equation}
We have \textup{(i)}
$\operatorname{tr}\{(\Eb^*_{b,K})^{-1}\} = \operatorname{tr}\{\Eb^{-1}_{b,K-2}\} + f_{b,i,j}(c_b)$ and
\textup{(ii)} $f_{b,i,j}(c_b) < 2(N-3)^{-1}$.
\end{lemma}

\begin{pf}
We first replace $\Eb$, $\Eb_{11}$ and $\Eb_{22}$ in Lemma~\ref
{lemmaS1} by
$\Eb^*_{b,K}$,
\[
\Eb_{11} = \left[ \matrix{ N_b &
c_b
\vspace*{2pt}\cr
c_b & N_b }
 \right],
\]
and $\Eb_{b, K-2}$, respectively. This allows to verify (i). In
addition, we have from Lemma~\ref{lemmaS2} that
\begin{eqnarray*}
f_{b,i,j}(c_b) &= &\operatorname{tr}\bigl\{\bigl(\Eb^*_{b,K}
\bigr)^{-1}\bigr\} - \operatorname{tr}\bigl\{\Eb^{-1}_{b,K-2}\bigr\}
\leq \operatorname{tr}\bigl\{\bigl(\Eb^*_{b,K}\bigr)^{-1}\bigr\} - \operatorname{tr}\bigl\{
\bigl(\Eb^*_{b,K-2}\bigr)^{-1}\bigr\}
\\
&\leq& \operatorname{tr}\bigl\{\bigl(\Eb^*_{b,K}\bigr)^{-1}\bigr\} - \operatorname{tr}\bigl
\{\bigl(\Eb^*_{b,K-1}\bigr)^{-1}\bigr\} + \operatorname{tr}\bigl\{ \bigl(\Eb
^*_{b,K-1}\bigr)^{-1}\bigr\} - \operatorname{tr}\bigl\{\bigl(
\Eb^*_{b,K-2}\bigr)^{-1}\bigr\}
\\
& <& 2(N-3).
\end{eqnarray*}
This proves (ii).
\end{pf}

We now are ready to prove that $\Eb^*_K$ is similar to a block matrix.
This is done by considering the expression of $\Eb^*_{b,K}$ in (\ref
{Eq:Estar}).
We then will show that if $|c| > 3$, then $f_{b,i,j}(c_b) \geq
2(N-1)^{-1}$, which contradicted with Lemma~\ref{lemmaS3}(ii), a
necessary condition of Lemma~\ref{lemmaS2}. Since the same argument can
be applied after
permuting the rows and columns of $\Eb^*_{b,K}$, $\Eb^*_K$ must have
off-diagonal elements equal to $-1$ or $3$, and
is thus similar to a block matrix. We begin this procedure by deriving
some useful results.
With Lemma~\ref{lemmaS3}(i) and equation (\ref{Eq:Ej}), we have
\begin{eqnarray*}
\operatorname{tr}\bigl\{\bigl(\Eb^*_{b,K}\bigr)^{-1}\bigr\} &= &\operatorname{tr}\bigl\{
\Eb^{-1}_{b,K-2}\bigr\} + f_{b,i,j}(c_b)
\\
&\leq&\min\bigl\{\operatorname{tr}\bigl\{\Eb^{-1}_{b,K-2}\bigr\} +
f_{b,i,i}(3_b), \operatorname{tr}\bigl\{\Eb ^{-1}_{b,K-2}
\bigr\} + f_{b,j,j}(3_b) \bigr\}.
\end{eqnarray*}
Thus, $f_{b,i,j}(c_b) \leq\min\{f_{b,i,i}(3_b), f_{b,j,j}(3_b)\}$.
Let $z_{b,g}(c_b) = \sqrt{z_{b,i,i}(c_b)z_{b,j,j}(c_b)}$; then
\begin{eqnarray*}
&&\frac{z_{b,i,i}(3_b) b_{b,j,j} + z_{b,j,j}(3_b)b_{b,i,i}}{2}\\
&&\qquad\geq \sqrt{z_{b,i,i}(3_b)
b_{b,j,j}z_{b,j,j}(3_b)b_{b,i,i}}
\\
&&\qquad= \sqrt{z_{b,i,i}(3_b)z_{b,j,j}(3_b)}
\sqrt{b_{b,i,i}b_{b,j,j}} \geq z_{b,g}(3_b)|b_{b,i,j}|.
\end{eqnarray*}
The last inequality is due to the Cauchy--Schwarz inequality
[see also, Theorem~14.10.1 of \citet{Harville1997bk31}]. With the same
reason, we also have $a^2_{b,i,j} \leq a_{b,i,i}a_{b,j,j}$. Thus,
for $|c| > 3$,
\begin{eqnarray*}
\bigl|z_{b,i,j}(c_b)\bigr| &=& |c_b - a_{b,i,j}|
\geq|c| - \bigl|N^{-1}\bigr| - |a_{b,i,j}| > 3_b -
\sqrt{a_{b,i,i}a_{b,j,j}}
\\
&\geq&3_b - \frac{a_{b,i,i} + a_{b,j,j}}{2} = \frac{3_b - a_{b,i,i} +
3_b - a_{b,j,j}}{2}\\
& =&
\frac{z_{b,i,i}(3_b) + z_{b,j,j}(3_b)}{2}.
\end{eqnarray*}
We note that $a_{b,i,i}$, and $a_{b,j,j}$ are positive since $\Eb
_{K-2}$ is positive definite.
Let $z_{b,a}(c_b) = (z_{b,i,i}(c_b) + z_{b,j,j}(c_b))/2$. We have, for
$|c| > 3$,
$|z_{b,i,j}(c_b)| > z_{b,a}(3_b)\geq z_{b,g}(3_b)$. It also can be
easily seen that
\[
f_{b,j,j}(3_b) = \frac{2[A_{b,j,j} + (N-3)b_{b,j,j}]}{(N-3)[A_{b,j,j} +
z_{b,j,j}(3_b)]}.
\]
With these facts and some algebra similar to that in \citet
{SatheShenoy1989ar2021},
we can see that, for $|c| > 3$,
\begin{eqnarray*}
&&\bigl[A_{b,i,i}A_{b,j,j} - z^2_{b,i,j}(c_b)
\bigr] f_{b,i,j}(c_b)
\\
&&\qquad= \bigl[(N-3)^2 +2(N-3)z_{b,a}(3_b)+
z^2_{b,g}(3_b) - z^2_{b,i,j}(c_n)
\bigr]f_{b,i,j}(c_b)
\\
&&\qquad\geq2\bigl[N-3 + z_{b,a}(3_b)\bigr]\\
&&\qquad\quad{} +
\bigl[(N-3)z_{b,a}(3_b) + z^2_{b,g}(3_b)
- z^2_{b,i,j}(c_b)\bigr] f_{b,i,j}(c_b).
\end{eqnarray*}
The first equality is due to $A_{b,i,j} = N - c + z_{b,i,j}(c_b)$. This
in turn leads to
$f_{b,i,j}(c_b) \geq2(N-3)^{-1}$. With Lemma~\ref{lemmaS3}(ii), we
thus can conclude that $|c| \leq3$ and
$\Eb^*_K$ is similar to a block matrix. The proof of Lemma~\ref
{lemmaA4} is then completed by using Lemma~\ref{lemmaA3}.
\end{appendix}
%\bibliographystyle{imsart-nameyear}

% imsref loaded by akundreckaite, 2015-07-09 12:35:54
%

%
%\begin{appendix}
%\section{}
%\end{appendix}

% zodis "Acknowledgments" paliekamas pagal autoriu
%\section*{Acknowledgments}

%\begin{supplement}[id=suppA]
%\sname{Supplement A}
%\stitle{}
%\slink[doi]{10.1214/00-AOSXXXXSUPP} %[doi,text={...}] - jei reikia
%suskaldyti doi
%\sdatatype{.pdf}
%\sfilename{aosXXXX\_supp.pdf}
%\sdescription{}
%\end{supplement}

%\begin{thebibliography}{99}
%\bibitem[\protect\citeauthoryear{}{}]{r1}
%\bibitem{r1}
%\end{thebibliography}

\printaddresses
\end{document}